\documentclass[a4paper]{article}
\usepackage{natbib}

  \usepackage{graphics,graphicx}
   \usepackage{amsmath}
   \usepackage{amsfonts}
   \usepackage{float,txfonts,amssymb,bm}
 \usepackage{colortbl}
\newcommand{\orderNcube}{\ensuremath{ O \left( \frac{1}{N^3}\right)}}

 \bmdefine\bone{1}    \bmdefine\bX{X}  \bmdefine\bt{t}
\bmdefine\btau{\tau} \bmdefine\bx{x}  \bmdefine\balpha{\alpha}
\bmdefine\bbeta{\beta} \bmdefine\by{y} \bmdefine\bu{u}
\bmdefine\bG{G} \bmdefine\bF{F}

\begin{document}

\title{Elemental estimators for the Generalized Extreme Value tail}
\author{Allan McRobie \\
Cambridge University Engineering Department\\
Trumpington St, Cambridge, CB2 1PZ, UK \\
fam20@cam.ac.uk}

\maketitle

\begin{abstract}
In a companion paper (McRobie(2013) arxiv:1304.3918), a simple set of `elemental' estimators was presented for the
Generalized Pareto tail parameter. Each elemental estimator: involves only three log-spacings;
is absolutely unbiased for all values of the tail parameter; is location- and scale-invariant; and is valid for all sample sizes $N$, even as small as $N= 3$.
It was suggested that linear combinations of such elementals could then be used to construct efficient unbiased estimators.
In this paper, the analogous mathematical approach is taken to the Generalised Extreme Value (GEV) distribution.
The resulting elemental estimators, although not absolutely unbiased, are found to have very small bias, and may thus provide
a useful basis for the construction of efficient estimators.
\end{abstract}

\maketitle

\section{Introduction}
Together with the Generalized Pareto Distribution (GPD), the Generalized
Extreme Value (GEV) distribution is central in extreme
value theory. Each distribution has three parameters $(\mu, \sigma, \xi)$
corresponding to location, scale and tail (or shape) respectively, and the estimation
of the tail parameter is a problem that has received much attention.
\cite{McRobieGPD} presented an elegant set of `elemental' estimators for the GPD tail parameter.
These estimators: involved three log-spacings; are location- and scale-invariant; and are absolutely unbiased for all tail parameters $-\infty < \xi < \infty$ and all sample sizes $N \geq 3$.
The idea was that linear combinations of such elemental estimators could then be constructed which would be efficient, whilst preserving the desirable properties of
lack of bias and small sample validity. A variety of linear combinations were considered, and although consistency proofs have yet to be constructed,
numerical evidence was presented which suggested that, for distributions
within the GPD family, the root mean square error for at least one simple linear combination converged in proportion to $1/\sqrt{N}$ for all $\xi$.

The idea of this paper is to construct the equivalent elemental estimators for the GEV, using the same mathematical approach as was adopted for the GPD.
Here we find that the resulting elementals, rather than being absolutely unbiased, have a very small bias, typically an order of magnitude smaller than the root
mean square error. As such, they may thus provide a useful basis within which linear combinations may provide
efficient estimators.

In both the GPD and GEV cases, the ultimate aim is of course the application to the wider domain of attraction case, using
appropriately chosen maxima from any distribution, rather than from the pure GPD or GEV families.
That extension is left for later consideration. Here, only the results for the pure GEV case are presented.

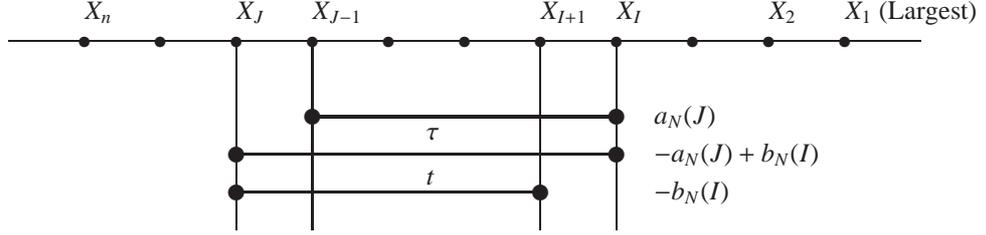
\begin{figure}[t!]
\begin{picture}(150,40)
\thinlines \multiput(30,30)(10,0){11}{\circle*{1.5}} \put (20,30)
{\line(1,0){120}} \put (30,33){$X_n$} \put(50,33){$X_J$}
\put(60,33){$X_{J-1}$} \put(90,33){$X_{I+1}$} \put(100,33){$X_I$}
\put (120,33){$X_2$} \put(130,33){$X_1$ (Largest)}
\put(50,5){\line(0,1){25}}\put(60,5){\line(0,1){25}}
\put(90,5){\line(0,1){25}}\put(100,5){\line(0,1){25}} \thicklines
\put(50,15){\line(1,0){50}}\thicklines\put(60,20){\line(1,0){40}}
\thicklines\put(50,10){\line(1,0){40}}
\put(50,15){\circle*{2}}\put(100,15){\circle*{2}}
\put(60,20){\circle*{2}}\put(100,20){\circle*{2}}
\put(50,10){\circle*{2}}\put(90,10){\circle*{2}} \put(105,19){$a_N(J)$}
\put(105,14){$-a_N(J)+b_N(I)$} \put(105,9){$-b_N(I)$} \put(75,11){$t$}
\put(75,17){$\tau$}
\end{picture}
\caption{The elemental construction for the GEV case is identical to that for the GPD case, except that different log-spacing weights
are used. } \label{bardiag}
\end{figure}

The elemental estimators are illustrated in Figure \ref{bardiag}. First, any sample from a GEV is ordered (with $X_1$ the data maximum and $X_N$ the data minimum). An elemental is then constructed from three log-spacings of the order statistics.
The construction takes any pair of non-adjacent order statistics $X_I$ and $X_J$, and involves the log-spacing
$\log( X_I- X_J)$  $(J \geq I+2)$ together with two smaller log-spacings $\log( X_{I+1}- X_J)$  and $\log( X_I- X_{J-1})$  that nestle within it.
Location- and scale-invariance is ensured by using only the statistics
\begin{equation}
\tau = \frac{X_I- X_{J-1}}{X_I-X_J} \ \ \  \text{and} \ \ \  t = \frac{X_{I+1}- X_J}{X_I-X_J} \label{taut}
\end{equation}
The resulting elemental estimator takes the form
\begin{equation}
\hat{\xi}_{IJN} = a_{N}(J) \log \tau   - b_{N}(I) \log t   \label{eqn1}
\end{equation}
where, as will be shown, the weights $a_{N}(J)$ and $b_{N}(I)$ are given by
\begin{eqnarray}
a_{N}(J) & = & - \left[ \left(
                  \begin{array}{c}
                    N \\
                    J-1 \\
                  \end{array}
                \right)
                \sum_{m=0}^{J-1}
 \left(
                  \begin{array}{c}
                    J-1 \\
                    m \\
                  \end{array}
                \right)
(-1)^m
\log(N-(J-1)+m)
 \right]^{-1} \label{anji} \\
 b_{N}(I) & = & - \left[ \left(
                  \begin{array}{c}
                    N \\
                    I \\
                  \end{array}
                \right)
                \sum_{m=0}^{I}
 \left(
                  \begin{array}{c}
                    I \\
                    m \\
                  \end{array}
                \right)
(-1)^m
\log(N-I+m)
 \right]^{-1}
 \label{ajbi}
\end{eqnarray}
These enjoy the property that $b_N(K-1) = a_N(K)$.
For the GPD case, the weights did not depend upon the sample size $N$, but for the GEV case, they do.

The Generalized Extreme Value Distribution arises as the limiting
distribution of block maxima (see for
example \citet{embrechts}). It has distribution function:
\begin{equation}
F(x) = \left\{  \begin{array}{ll}
 \exp \left( - \left( 1+ \xi z \right)^{-1/\xi} \right) & \ \ \text{for} \ \ \xi \neq 0 \\
 \exp \left( - \exp \left( -  z  \right) \right) & \ \ \text{for} \ \ \xi = 0
\end{array} \right. \ \ \text{where} \ \ z = \frac{x-\mu}{\sigma}
\label{thegpd}
\end{equation}
The parameters $\mu$ and $\xi$ can take any value on the real line,
whilst $\sigma$ can be any positive value.
For GEVs with positive $\xi$, the support ($\mu-\sigma/\xi \leq x$) is bounded below
but unbounded at the right.
For $\xi$ negative, the support is unbounded below but bounded above ($  x \leq \mu - \sigma / \xi$).
Since the right tail is the tail of interest, positive $\xi$ corresponds to long- or heavy-tailed distributions which are unbounded above, and negative $\xi$ corresponds to ``bounded-above'' distributions which have a finite right end point.

\section{Elemental Estimators for the GEV}

The construction of the elemental estimators for the GEV follows that for the GPD case presented in \cite{McRobieGPD}.
The detailed construction is given in Appendix 1 here. In outline, a sample of size $N$ is drawn from a member of the pure GEV family.
The sample is then ordered, with $X_1$ the sample maximum and $X_N$ the sample minimum.
The location- and scale-invariant statistics $\tau$ and $t$  are then constructed, these being ratios of data spacings as defined in Fig.~\ref{bardiag} and
Eqn~\ref{taut}.

The expected values of $\log \tau$ and $\log t$ are then determined for the two cases where $\xi$ is positive or negative.
These integrals are expressed via the Probability Integral Transform as integrals of the uniformly-distributed distribution
function $F(\bX)$ over the $N$-dimensional unit simplex corresponding to the ordered sample.

These integrals readily decompose into a simple part (of the form  $\xi \langle \log(-\log F) \rangle$) and a complicated part (of the form
$ \langle \log(1-\phi^{|\xi|}) \rangle$). The simple part leads to a term  proportional to $\xi$, and forms the core of the estimator.
For the GPD, the elemental combination of $\langle \log \tau \rangle$ and $ \langle \log t \rangle$ whose simple terms deliver the estimate $\xi$
has the pleasing property that the complicated terms cancel exactly at all $\xi$ (such that the complicated integrals do not need to be evaluated explicitly).
However, in the GEV case, the elemental combination whose simple parts deliver the estimate $\xi$ does not cancel the complicated parts.
However, as will be demonstrated numerically, the residual error in this non-cancellation is
typically much smaller than the variability.

\begin{figure}[h!] \centering
  \includegraphics[width=70mm,keepaspectratio]{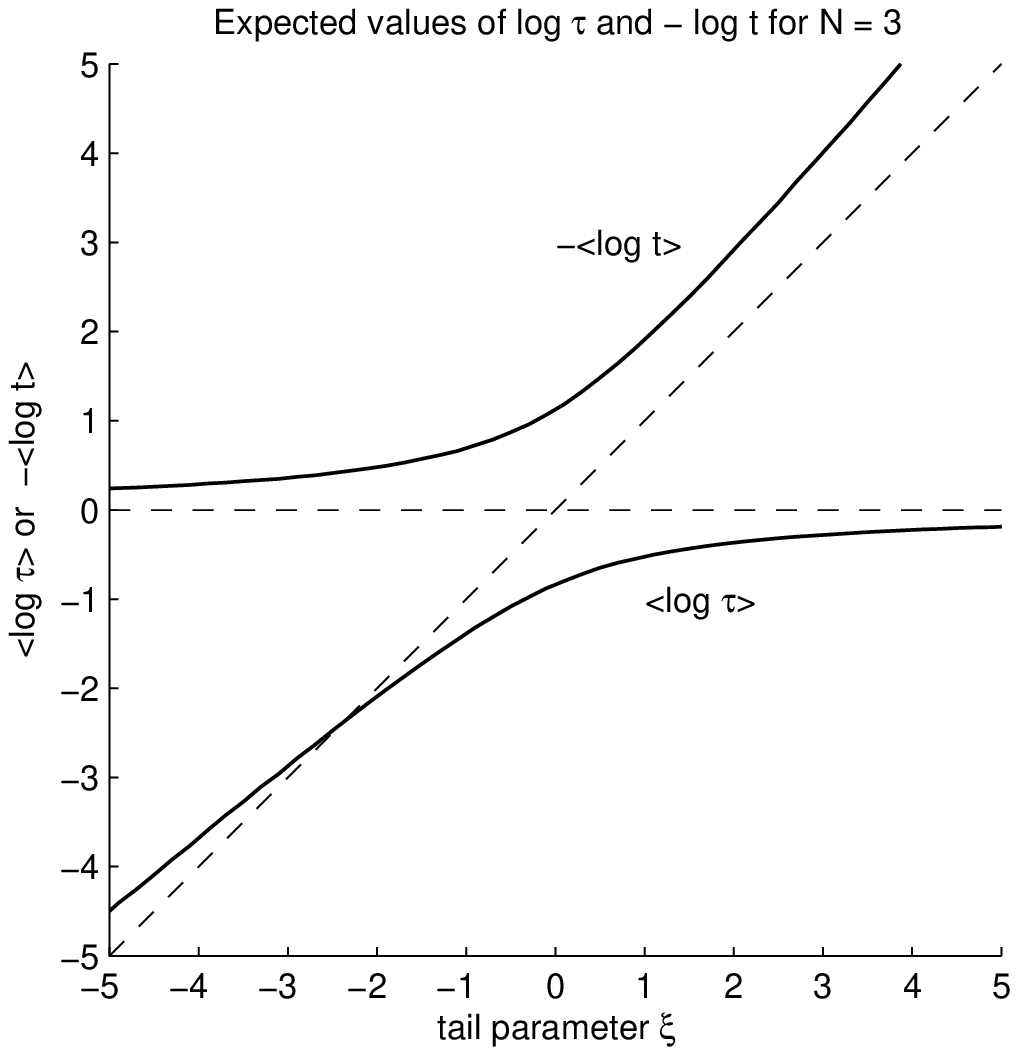}
  \includegraphics[width=70mm,keepaspectratio]{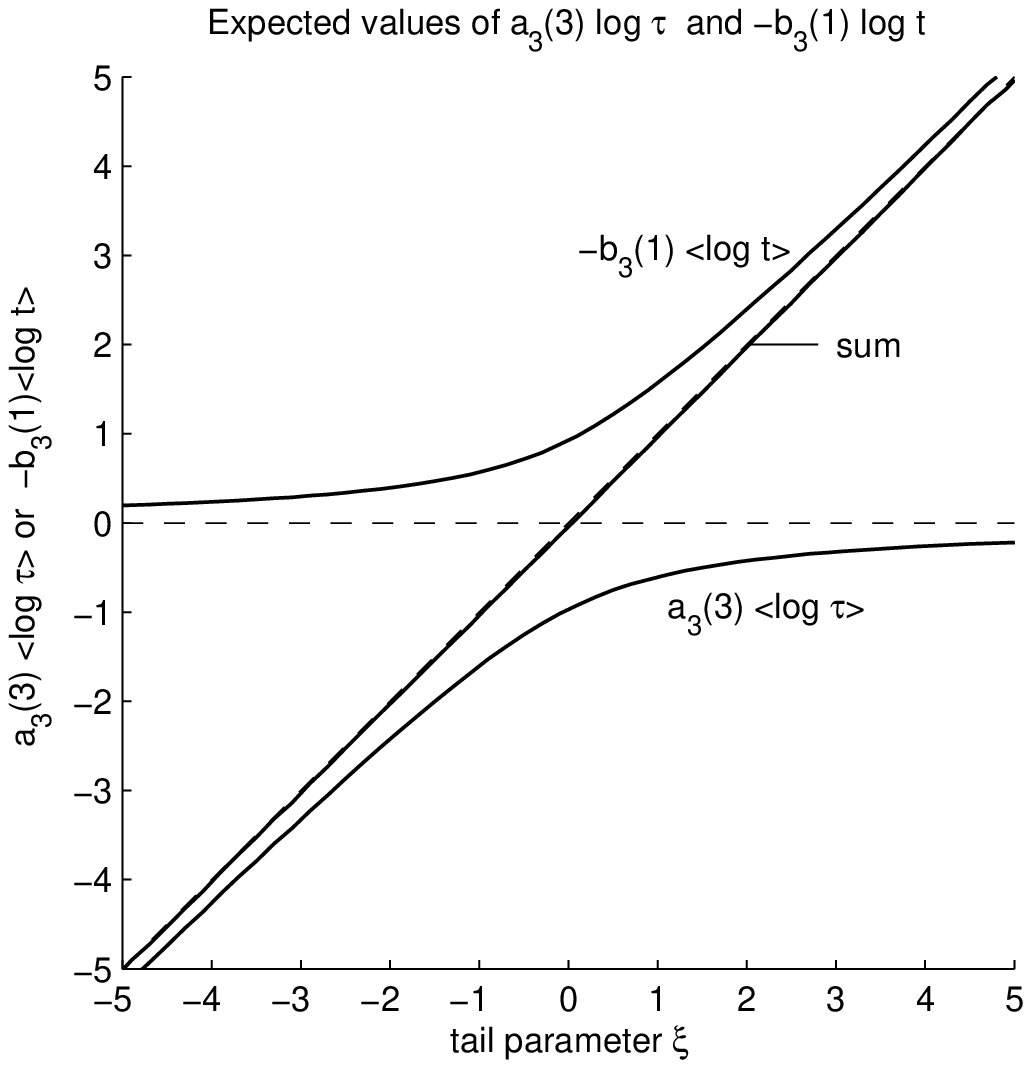} //
  \includegraphics[width=70mm,keepaspectratio]{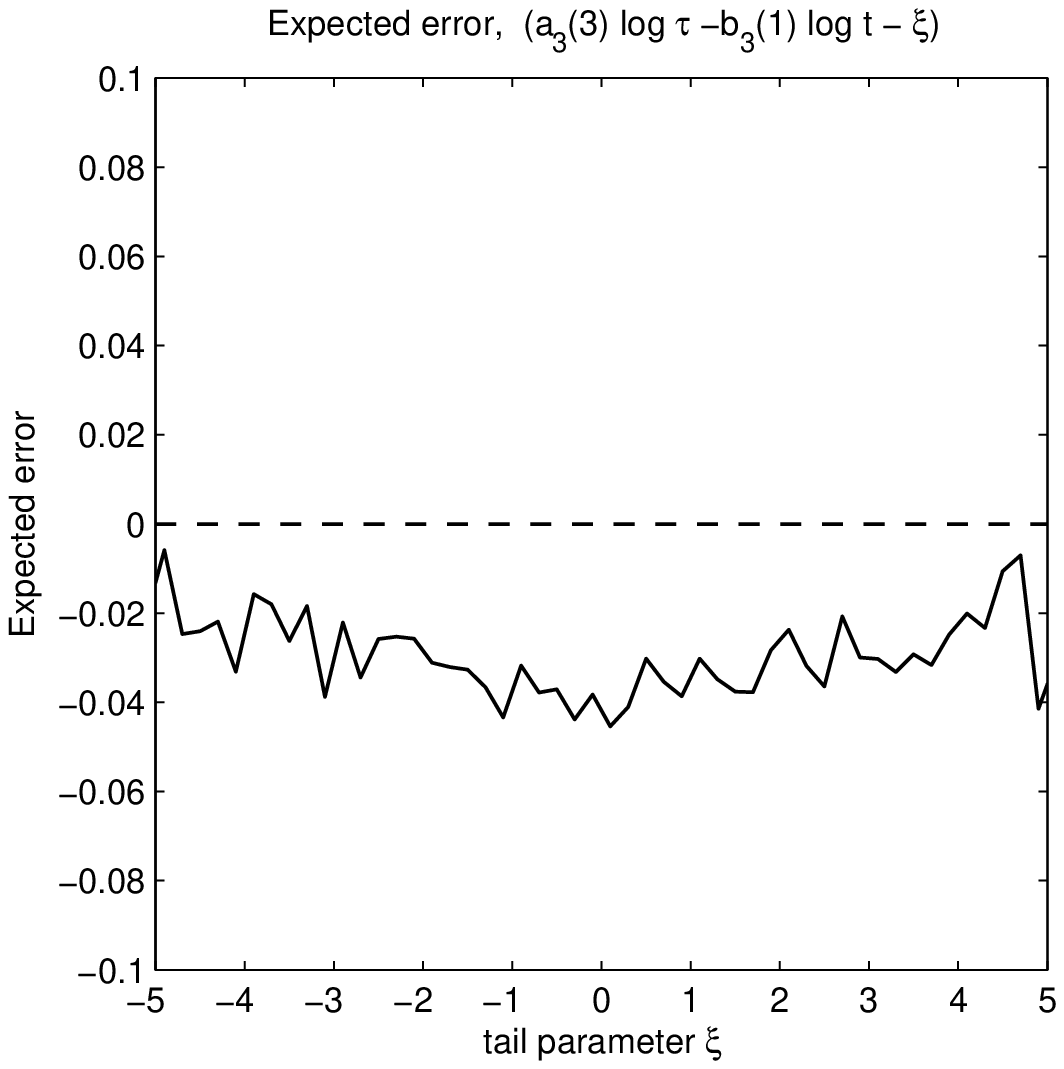}
  \includegraphics[width=70mm,keepaspectratio]{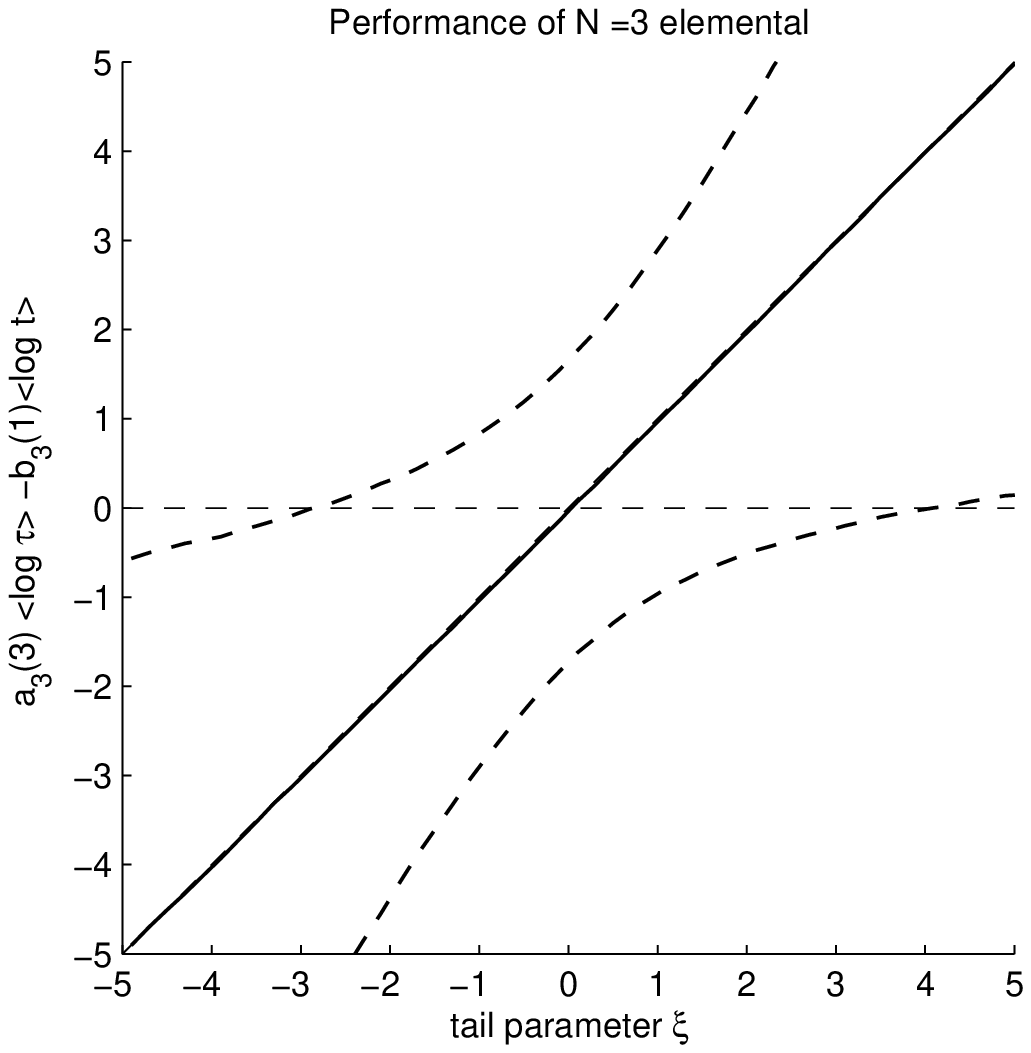}
 \caption{The elemental construction for $N=3$. The first figure shows how the expected values of $\log \tau$ and $-\log t$ have approximately hyperbolic form,
 and the second figure shows how scaling these by $a_{N}(J)$ and $b_{N}(I)$ makes each asymptotic to the diagonal for large $|\xi|$. Their sum then lies very close to the diagonal, but not exactly.
The third figure shows the small bias, and the fourth shows the standard deviation of the estimates.}\label{GEV2a}
\end{figure}

Although the construction thus far has been developed by simply mimicking the procedure
adopted for the GPD elementals, Fig.~\ref{GEV2a} shows how the construction is actually a natural one for the GEV per se.
For $N=3$, the first figure shows how the expected values of $\log \tau$ and $-\log t$ have approximately hyperbolic form,
asymptoting to zero at one end and to some linear slope at the other.
The second figure shows how scaling these by $a_N(J) = a_3(3) = -1/(3\log(3/4))$ and $b_N(I) = b_3(1) = -1/(3\log(2/3)) $ makes each asymptotic to the diagonal for large $|\xi|$. The elemental estimator is the
summation, and this lies very close to the diagonal as desired, but not exactly so.
It is asymptotic to the diagonal for large $|\xi|$ (as may be expected
since the expected values of the complicated terms $\log (1-\phi^{|\xi|})$ tend to $\log 1 = 0$ for $|\xi|$ large), and the small but non-zero bias
for intermediate values of $|\xi|$ is shown in the third figure.
The standard deviations of the estimator are shown in the fourth figure, from whence it can be seen that, at worst (near $\xi = 0$),
the bias is some 50 times smaller than the standard deviation.
These figures were produced numerically, using 250,000 samples of size $N=3$ at each value of $\xi$. Although
the figures correspond to the single elemental available at $N= 3$, similar figures apply for any elemental at any $N$.

In summary, the elementals for the GEV are constructed by choosing weights $a_N(J)$ and $b_N(I)$ which guarantee lack of bias at large $|\xi|$.
As stated, one hopes that the bias at any intermediate value of $\xi$ is small, and happily it is typically much smaller than the variability there
(and even if larger intermediate deviations were to occur for larger $N$, it is clear that linear combinations of elementals could be constructed
such as to manage the overall deviation).

%
%
%
%
%

%

\section{Evaluation of the Coefficients}
\begin{figure}[h!] \centering
  \includegraphics[width=70mm,keepaspectratio]{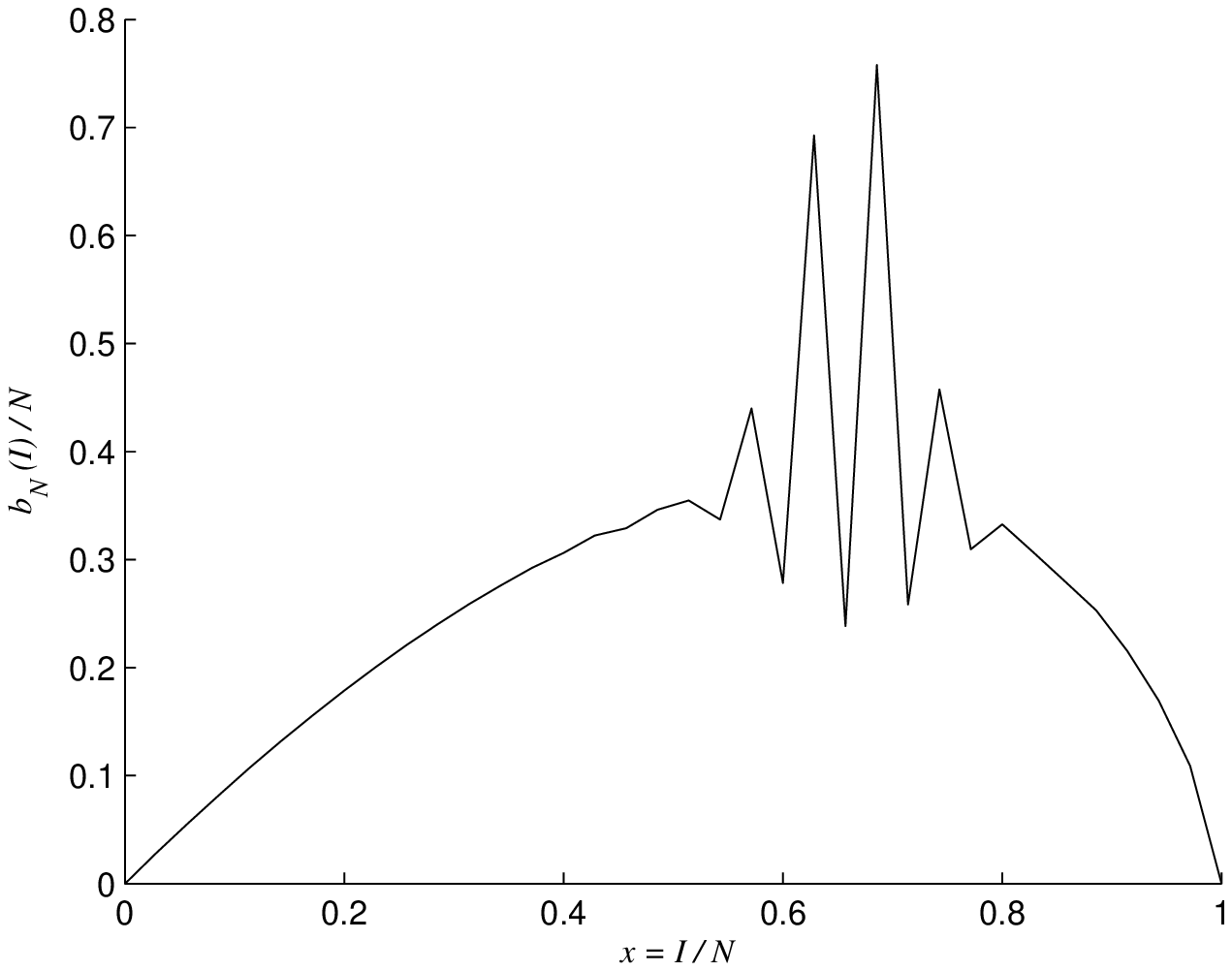}
  \includegraphics[width=70mm,keepaspectratio]{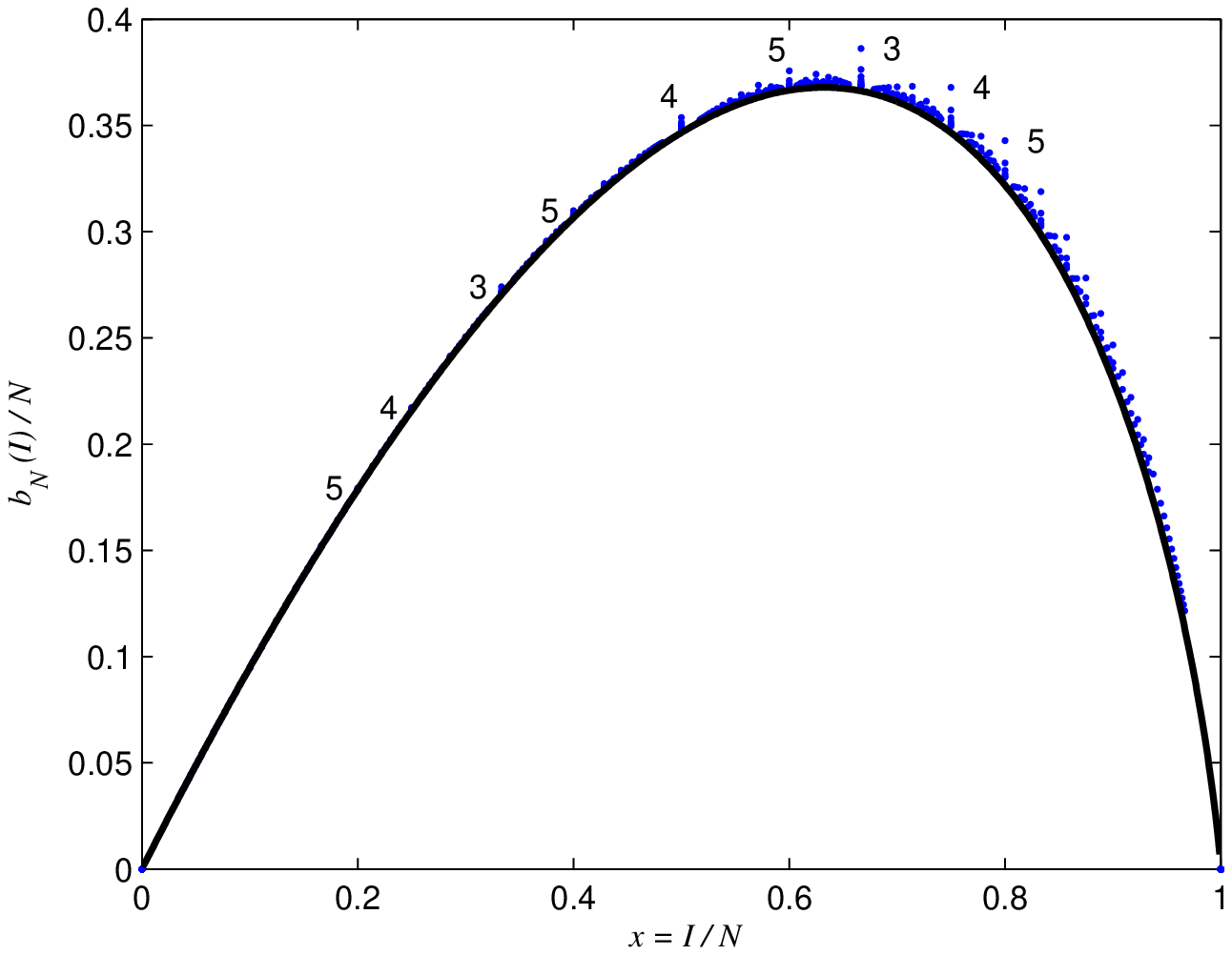}
  \caption{Calculation of the elemental coefficients $b_N(I)$. The left hand diagram illustrates erroneous numerical oscillations from the recursion formula, clearly visible for $N = 35$.
  The right hand diagram shows the convergence of $b_N(I)/N$ (shown as dots for $N= 3$ to 30) to $f(x) = -(1-x)\log(1-x)$ (shown solid), where $x = I/N$.
  }\label{Acalc}
\end{figure}

Although the elemental coefficients $a_N(J)$ and $b_N(I)$ are given by a comparatively simple summation
(Eqns.~\ref{ajbi} and \ref{ajbi}) involving binomial coefficients
and logarithms, the evaluation becomes fraught with numerical errors for $N \gtrsim 25$ in 32 bit precision computation.
The numerical oscillations are illustrated in Fig.~\ref{Acalc}. That this is a numerical instability rather than true behaviour is noted
from the fact that if one makes less effort to calculate the factorials carefully then the instability occurs at lower $N$.
Even when care is taken, and the binomial coefficients are calculated via the standard procedure of using logarithms of gamma functions, there is only marginal improvement  in
numerical performance. The problem arises because the summations in Eqns.~\ref{anji} and \ref{ajbi} sum close to zero, but the individual elements in
the sums - the logarithms weighted by the binomial coefficients - may be large. Repeated addition and subtraction of these,
via the $(-1)^m$ factor, means that small rounding errors accumulate and eventually dominate.

Before describing how these numerical instabilities can be resolved, we note that the coefficients $a_N(J)$ follow trivially from the
$b_N(I)$ coefficients using the relation $a_N(J) = b_N(J-1)$  (with the nuance that
$b_N(I)$  must also be evaluated at the meaningless value $I = N-1$ in order
to capture the coefficient $a_N(N)$ which would not otherwise be computed).
We thus proceed to describe only the calculation of the $b_N(I)$.

We define $\beta_N(I) = 1/b_N(I)$, such that $\beta$ is the weighted summation and its reciprocal $b$ is the desired coefficient.
It follows from Eqn.~\ref{ajbi} that the coefficients $\beta$ can be computed using the recursion relation
\begin{equation}
\beta_N(I+1) = \frac{N}{I+1} \beta_{N-1}(I) -\frac{N-I}{I+1}\beta_{N}(I)
\label{recursion}
\end{equation}
This is seeded by the $I = 1$ top row
\begin{equation}
\beta_N (1)  = N \log \left( 1- \frac{1}{N} \right)
\end{equation}
Each $\beta$ coefficient in row $I+1$ is thus computed from a weighted sum of
two terms above it in row $I$, much akin to a standard Pascal's Triangle construction.
This elegant construction
avoids all need for explicit computation of factorials or logarithms of gamma functions.
However, even this algorithm does not remove the numerical inaccuracies that
arise for $ N \gtrsim 25$, thus an additional approach is required for large $N$.

It can be readily seen that, to first order (ignoring terms of order $1/N$), the recursion formula reduces to
\begin{equation}
\beta_{N-1}(I) =  (1-x) \beta_N(I) + x \beta_{N}(I+1)
\label{lininterp}
\end{equation}
where $ x = I/N$. This is simple linear interpolation amongst the three coefficients. 
To first order, then, the approximate recursion formula Eqn.~\ref{lininterp}
states that the three values of $\beta$ lie on a straight line
 when plotted as a function of $x$. This suggests that there is some underlying curve to which the $\beta$ values are converging.

Numerical evidence presented in Fig.~\ref{Acalc} b)  suggests that
\begin{equation}
\lim_{N \rightarrow \infty}  \ \ \ \left[   N \  \left(
                  \begin{array}{c}
                    N \\
                    I \\
                  \end{array}
                \right)
                 \sum\limits_{m=0}^{I}
 \left(
                  \begin{array}{c}
                    I \\
                    m \\
                  \end{array}
                \right)
(-1)^m
\log(N-I+m) \right]
  =   \frac{1}{(1-x)\log(1-x)} \ \ \ \ \text{where} \ \ \  x = \frac{I}{N}
  \label{curious}
\end{equation}
such that $b_N(I)/N$ converges to  $f(x) = -(1-x)\log(1-x)$ as $N$
increases. An analytical demonstration of this is presented in Appendix 2.

Further numerical analysis suggests that, for finite $N$, there is a better approximation
\begin{equation}
\frac{b_N(I)}{N}  \approx  - (1-x) \log(1-x) -  \frac{x}{12 N} \log(1-x)
\label{bNapprox2}
\end{equation}

\begin{figure}[h!] \centering
  \includegraphics[width=120mm,keepaspectratio]{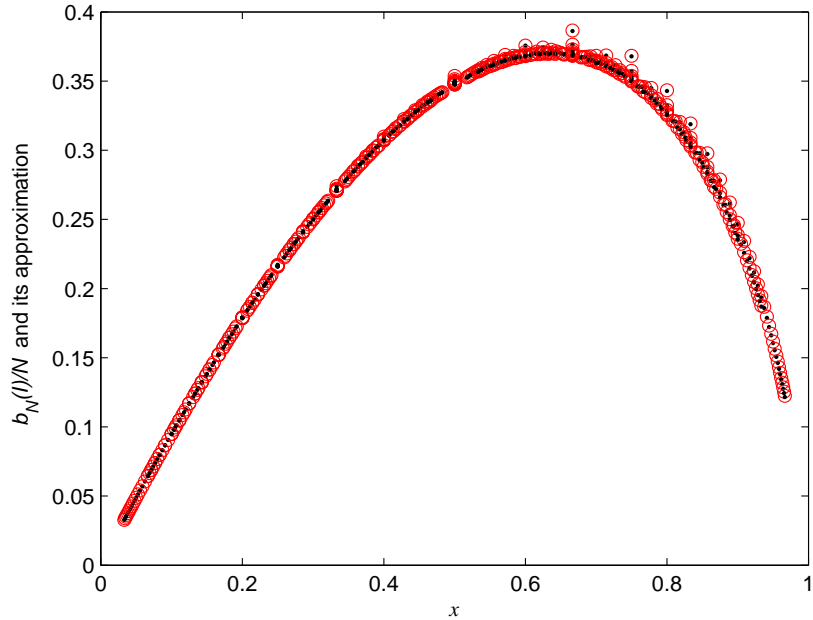}
  \caption{The approximation of $b_N(I)/N$ by Equation~\ref{bNapprox2}. The $b_N(I)/N$ values are plotted as points and the approximations are plotted using small circular markers.
  }\label{Acalc1a}
\end{figure}

\begin{figure}[h!] \centering
  \includegraphics[width=60mm,keepaspectratio]{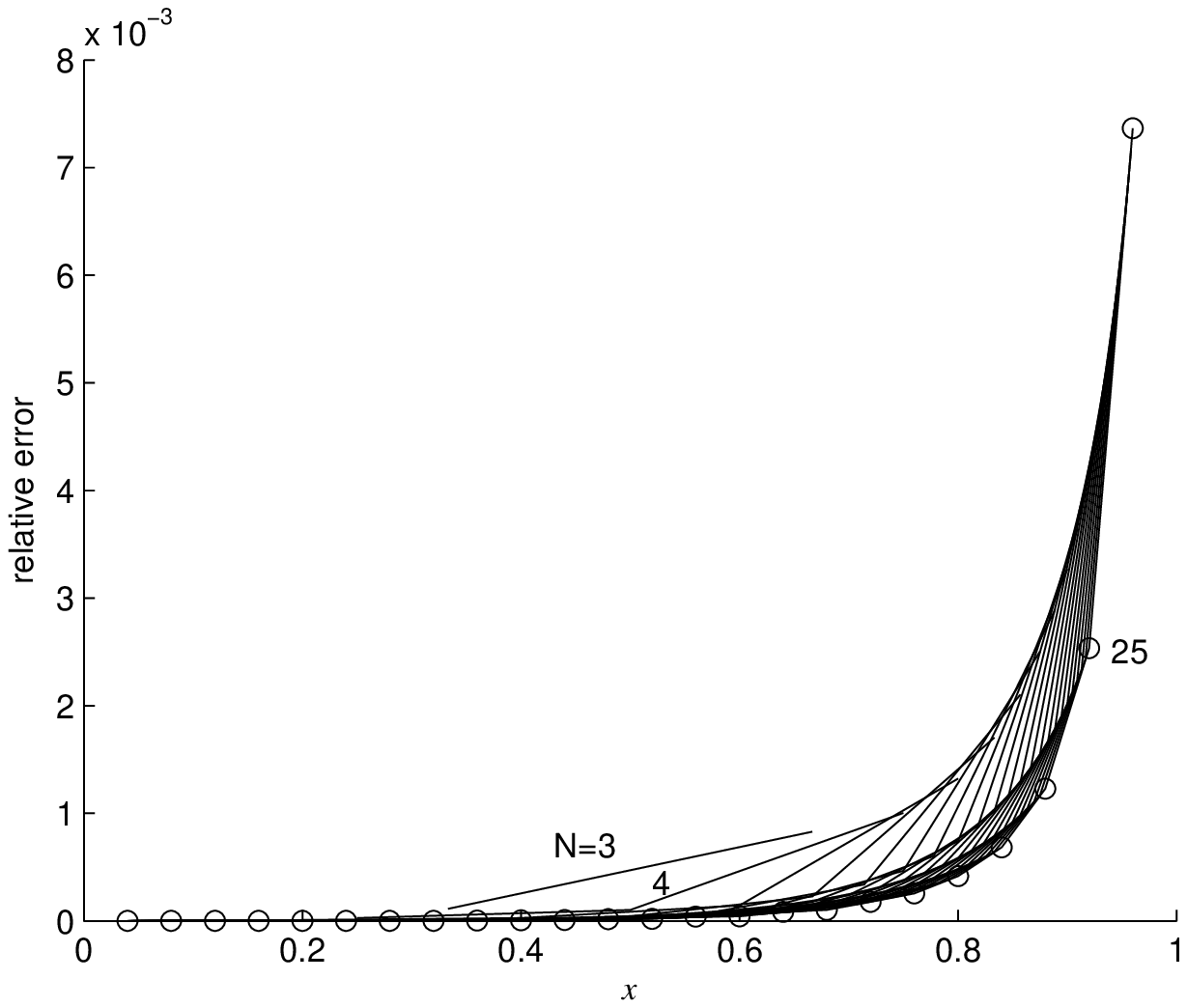}
  \includegraphics[width=60mm,keepaspectratio]{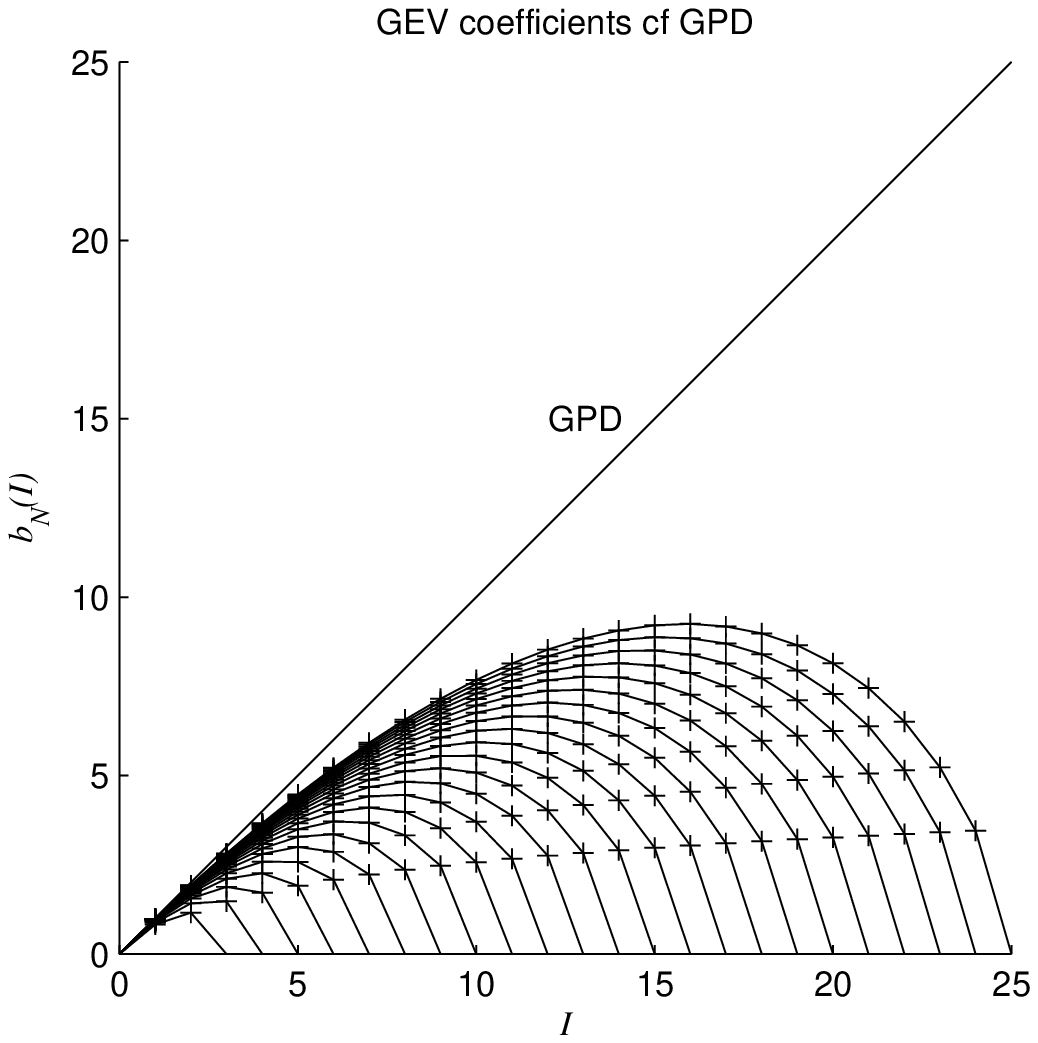}
  \caption{Left: the relative error involved in the approximation of $b_N(I)/N$ by Equation~\ref{bNapprox2}. Right: a comparison of
  the elemental coefficients $b_N(I)$ of the GEV with those of the GPD.
  }\label{Acalc1b}
\end{figure}

The numerical evidence for this is shown in Fig.~\ref{Acalc1a} where, for $N$ up to 30, the actual values of $b_N(I)/N$ are plotted as points and
the approximations are plotted using small circular markers. In all cases, even for $N=3$, each actual value lies well within the corresponding circular marker.
The relative error $(b_N^{approx}-b_N)/b_N$ is shown in Fig.~\ref{Acalc1b}, showing convergence at any $I$ and a better than 1 \% accuracy for all $N \leq 25$
 (and in most cases shown, the accuracy is considerably better).

In summary, for $N \leq 25$, the $b_N(I)$ are computed by the recursion formula Eqn.~\ref{recursion},
and for $N>25$ we use the approximation for $b_N(I)$ contained within Eqn.~\ref{bNapprox2}. The corresponding values of
$a_N(J)$ are computed from the relation that $a_N(K) = b_N(K-1)$.

For the GPD case \cite{McRobieGPD}, the elemental coefficients may be succinctly expressed as
\begin{equation}
a_N^{GPD}(J) = J-1  \ \ \ \text{and} \ \ \  b_N^{GPD}(I) = I
\end{equation}
Unlike the GEV case, these are independent of $N$.
Fig.~\ref{Acalc1b}b) shows that the GEV coefficients
converge towards the GPD weightings for the extreme tails (i.e. for $x$ small), much as may be expected.

\section{Performance of the elementals}

\begin{figure}[h!] \centering
  \includegraphics[width=70mm,keepaspectratio]{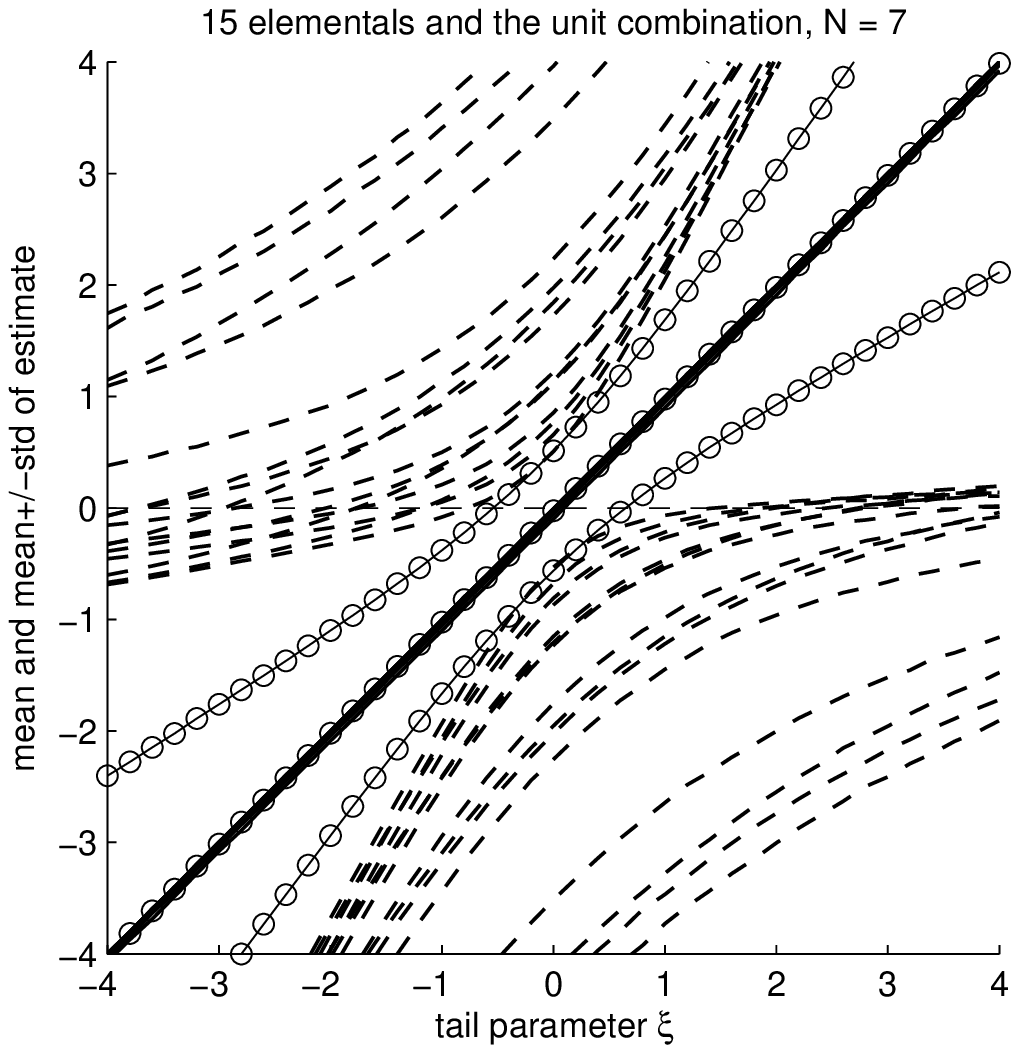}
  \includegraphics[width=70mm,keepaspectratio]{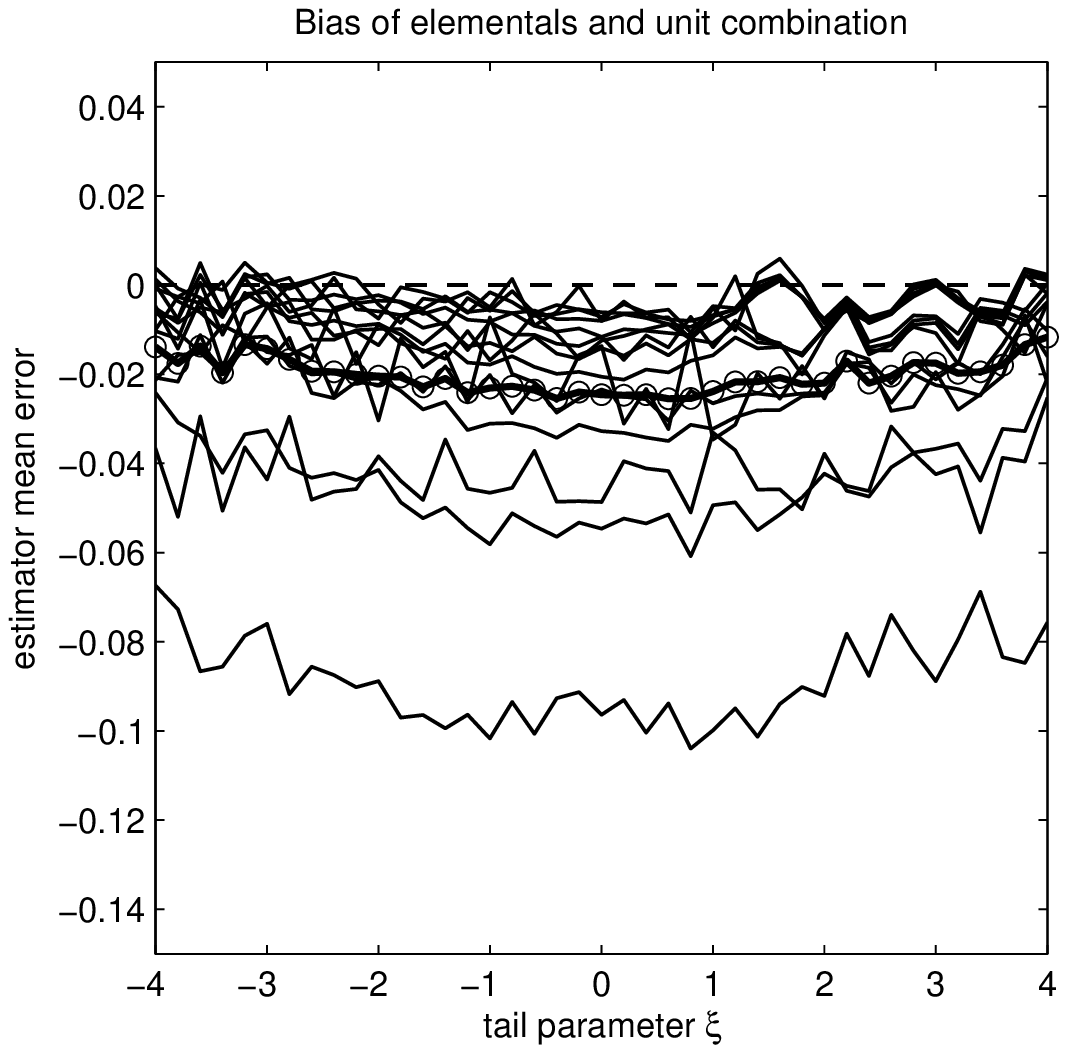}
  \caption{The left hand diagram shows the performance (mean, and mean $\pm$ standard deviation) of each of the fifteen elementals
  available for samples of size $N=7$. Also shown (circles) is the performance of the unit sum
  linear combination which gives equal weight to each elemental.
  The right hand diagram shows the small bias of the elementals and the linear combination. Plots based on 500,000 samples at each $\xi$.
  }\label{GEVN7afig1}
\end{figure}

Fig.~\ref{GEVN7afig1}a) shows the performance (mean and mean $\pm$ std. dev) of each of the
fifteen elemental estimators available at $N=7$ when applied to samples drawn from a pure GEV.
It can be seen that each elemental is very close to being unbiased. It can also be seen that some
standard deviations are large.
The biases are plotted in Fig.~\ref{GEVN7afig1}b).
Also plotted in each of these figures is the
performance of the unit sum linear combination which gives equal weight to each of the fifteen elementals. As may be expected, the linear combination - like the elementals - is almost unbiased, and generally has
a considerably smaller standard deviation than the individual elementals.

\begin{figure}[h!] \centering
  \includegraphics[width=70mm,keepaspectratio]{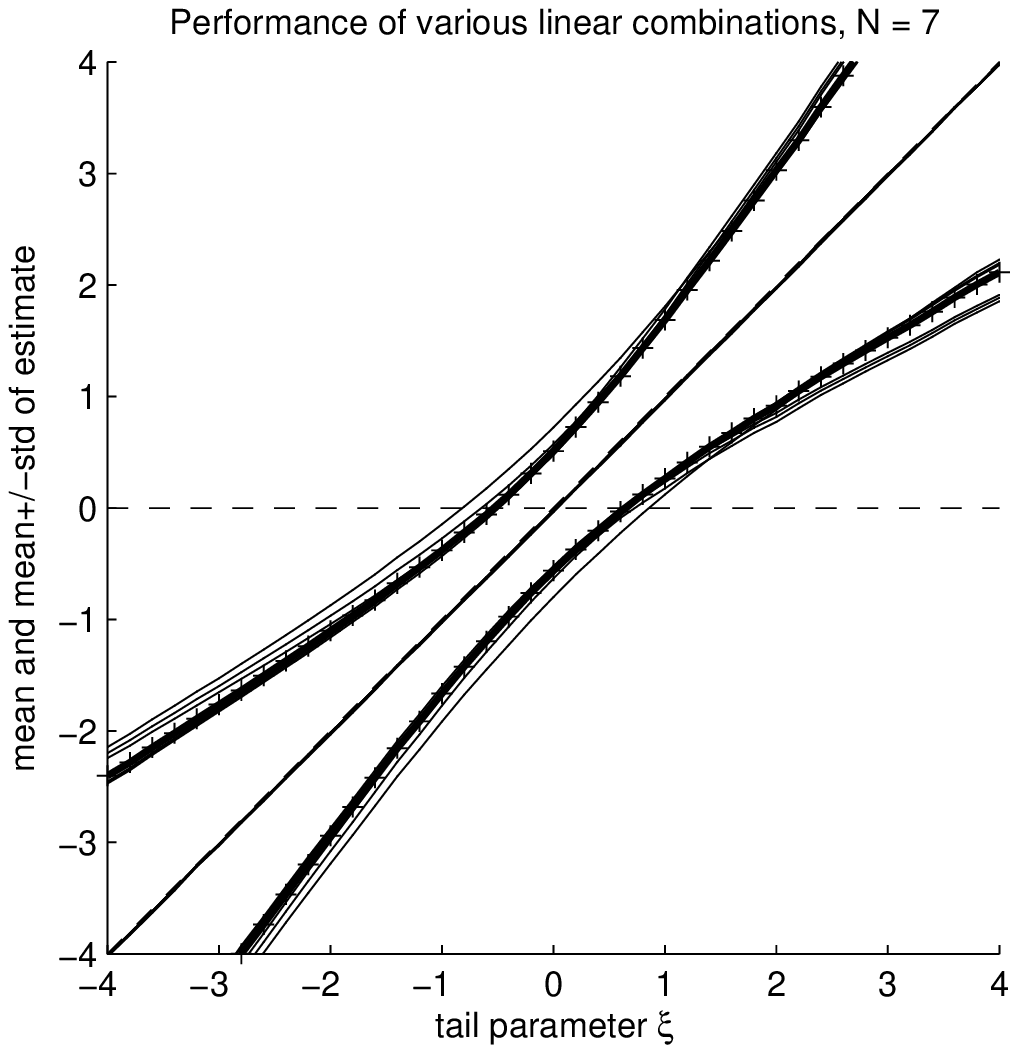}
  \includegraphics[width=70mm,keepaspectratio]{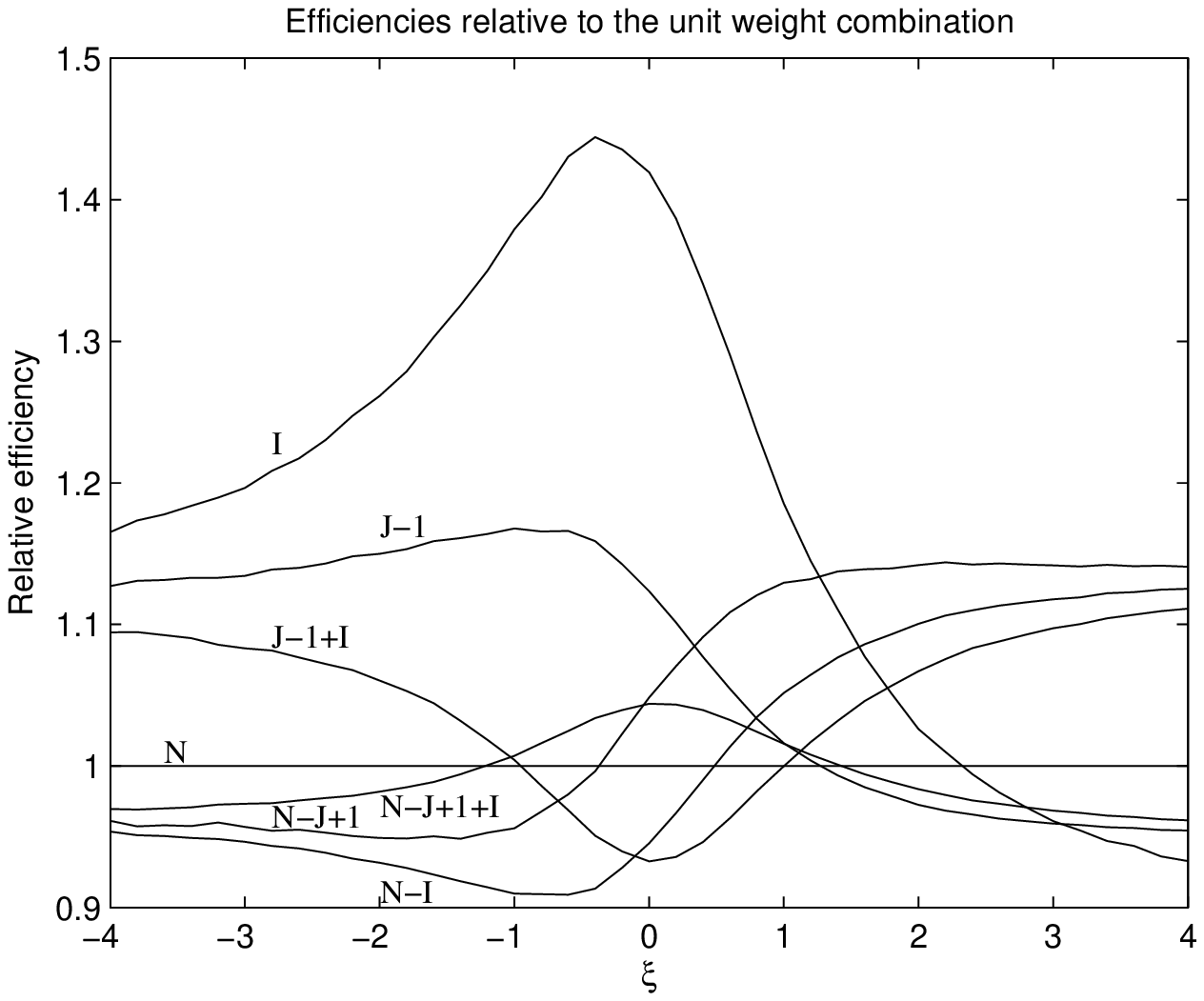}
  \caption{The left hand diagram shows the performance (mean (solid), and mean $\pm$ standard deviation (dashed)) of seven different linear combinations of the fifteen elementals
  available for samples of size $N=7$. The equal weight combination is marked with $+$ signs.
  The right hand diagram shows the relative efficiency (ratio of root mean square errors) of each combination to the equal weight combination (denoted $N$). Plots based on 100,000 samples at each $\xi$.
  }\label{GEVN7bfig1}
\end{figure}

The question of which linear combination is in some sense optimal is addressed in Fig.~\ref{GEVN7bfig1} where  seven different linear combinations
of the individual elementals are considered. The weights in each linear combination have unit sum.  The first combination gives equal weight to each elemental.
The next three combinations resemble linear triangular basis function which - prior to normalisation
- have unit weight
at one corner of the upper triangular matrix of possible $I$ and $J$ values, and zero
weight at the other two corners. The weights are thus proportional to $N-J+1$, $J-1-I$ and $I$ respectively. The final three combinations are sums of two of the previous three, such that - pre-normalisation -  they have unit weight at two corners and zero weight at the remaining corner. There is nothing special about this choice of combinations, other than that the choice serves as a useful starting point.

 Fig.~\ref{GEVN7bfig1}a) shows that the combinations perform similarly. By plotting the ratio of root mean square errors,  Fig.~\ref{GEVN7bfig1}b) makes clearer the small differences in performance showing the efficiency of each combination relative to the equal weight combination. It can be seen that some combinations perform better in some parameter regimes, and that no combination is optimal for all $\xi$.
Whilst one could consider constructing some form of adaptive estimator -  a first estimate of $\xi$ being used to select the linear combination that is efficient in that region -
at this stage of the analysis, we simply proceed with any convenient choice of combination, such as that with equal weights, or - to mirror the choice used in the GPD case - the combination whose weights are proportional to $N-J+1$.

\section{Consistency}
No attempt is made here at providing an analytical proof of consistency for any combination. One reason for this is that small and moderately-sized data sets are often of practical interest, and the usefulness of a consistency proof requiring very large sample sizes is then debatable.

\begin{figure}[h!] \centering
  \includegraphics[width=70mm,keepaspectratio]{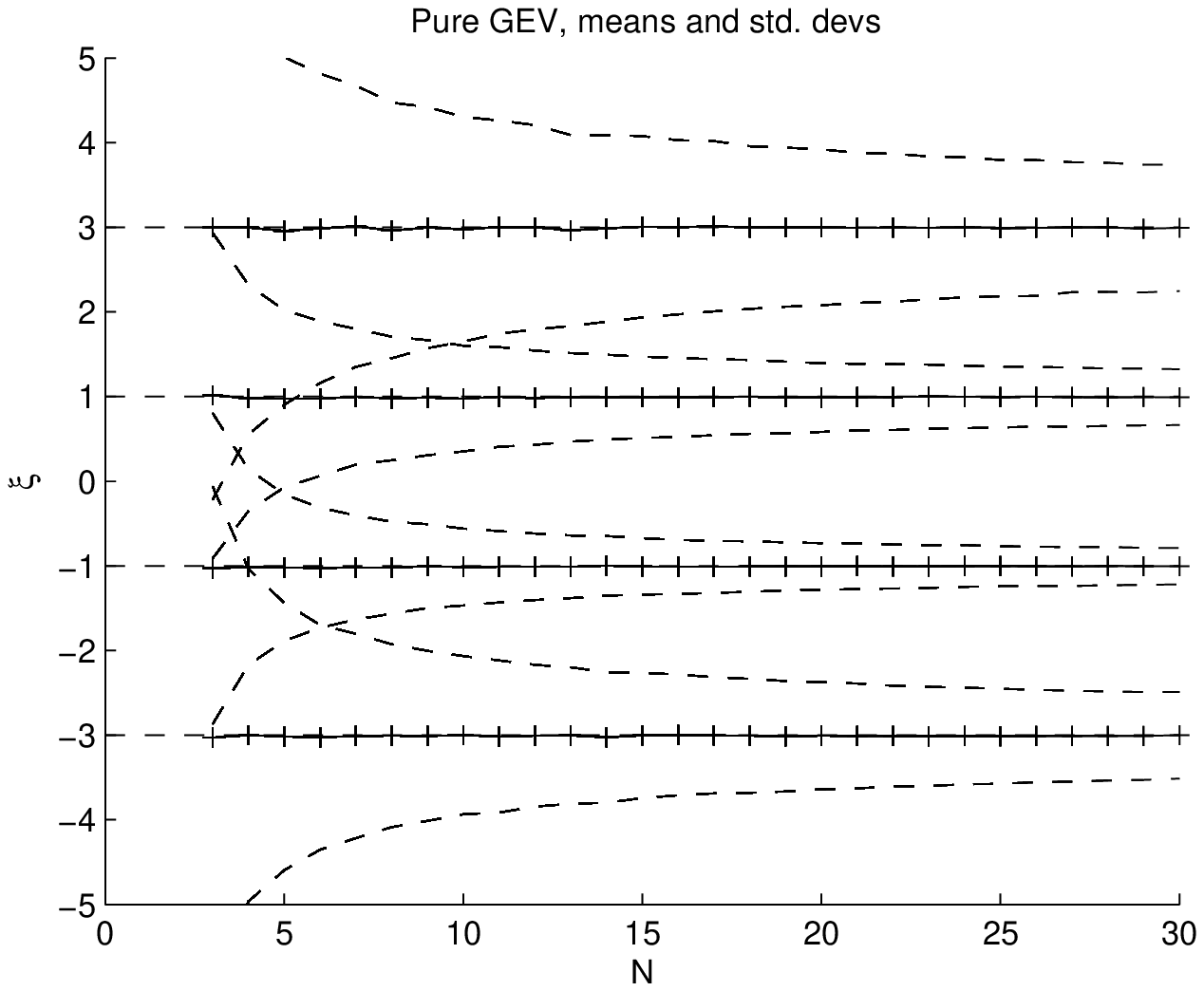}
  \includegraphics[width=70mm,keepaspectratio]{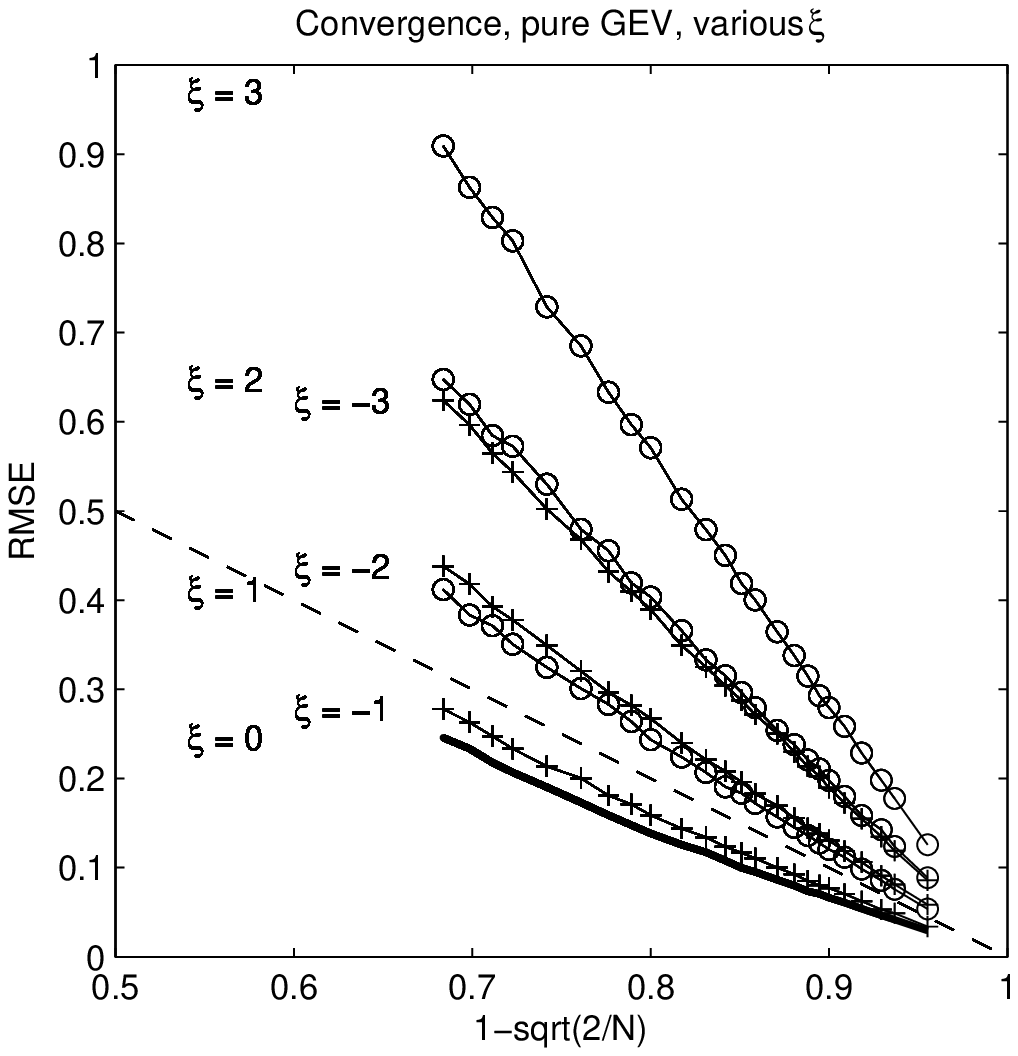}
  \caption{The left hand diagram shows the performance (mean, and mean $\pm$ standard deviation) of the elemental combination with weights
 proportional to $N-J+1$, for sample sizes from $N=3$ to $30$. Means are shown (+) with deviations dashed.
 The right hand diagram extends this to samples sizes up to $N=1000$. Cases with $\xi$ negative are shown (+) and those with
 $\xi$ positive are shown (o).
 The $y$-axis is the root mean square error (based on 10,000 samples) and the $x$-axis
 plots $1-\sqrt{2/N}$. The tendency towards linearity
 under increasing $N$ at each $\xi$
 suggests that there is consistency, with convergence at a rate proportional to $1/\sqrt{N}$.
 }\label{GEVconsistency}
\end{figure}

Instead, numerical evidence is provided to demonstrate that root mean square errors of some standard linear combinations appear to converge sensibly as sample sizes increase (with samples drawn from pure GEVs).
Fig.~\ref{GEVconsistency}a) shows how root mean square errors decrease at a selection of parameter values $\xi$ as sample size grows from $N=3$ to $30$, for the combination
whose weights are proportional to $N-J+1$.
This evidence is extended in Fig.~\ref{GEVconsistency}b) to samples of size $N=1000$. The apparent linearity of the graphs on the chosen axes
suggests that root mean square error decays as $1/\sqrt{N}$.

\section{Comparison with Maximum Likelihood}
Maximum Likelihood Estimation (MLE) is a standard tool for parameter estimation,
and \cite{embrechts} state that the numerical calculation of ML estimates for the
GEV ``poses no serious problem {\it in principle}''. Even with their italicised qualification, this is perhaps optimistic.
For example, it is well-known that MLE applied to the GEV requires regularity conditions that do not hold when $\xi <= -1/2 \ $ (\cite{smith}).
Moreover, for any data set, there are always parameter sets wherein the likelihood is arbitrarily large
(i.e.\ the likelihood is infinite everywhere on the boundary surface $1 +\xi (x_{max}-\mu)/\sigma = 0$ for $\xi <-1$, and is arbitrarily large in open
neighbourhoods adjacent to this).
\cite{castillo} also show that there are occasions where the likelihood function for the two-parameter GPD has no local maximum.

\begin{figure}[h!] \centering
  \includegraphics[width=70mm,keepaspectratio]{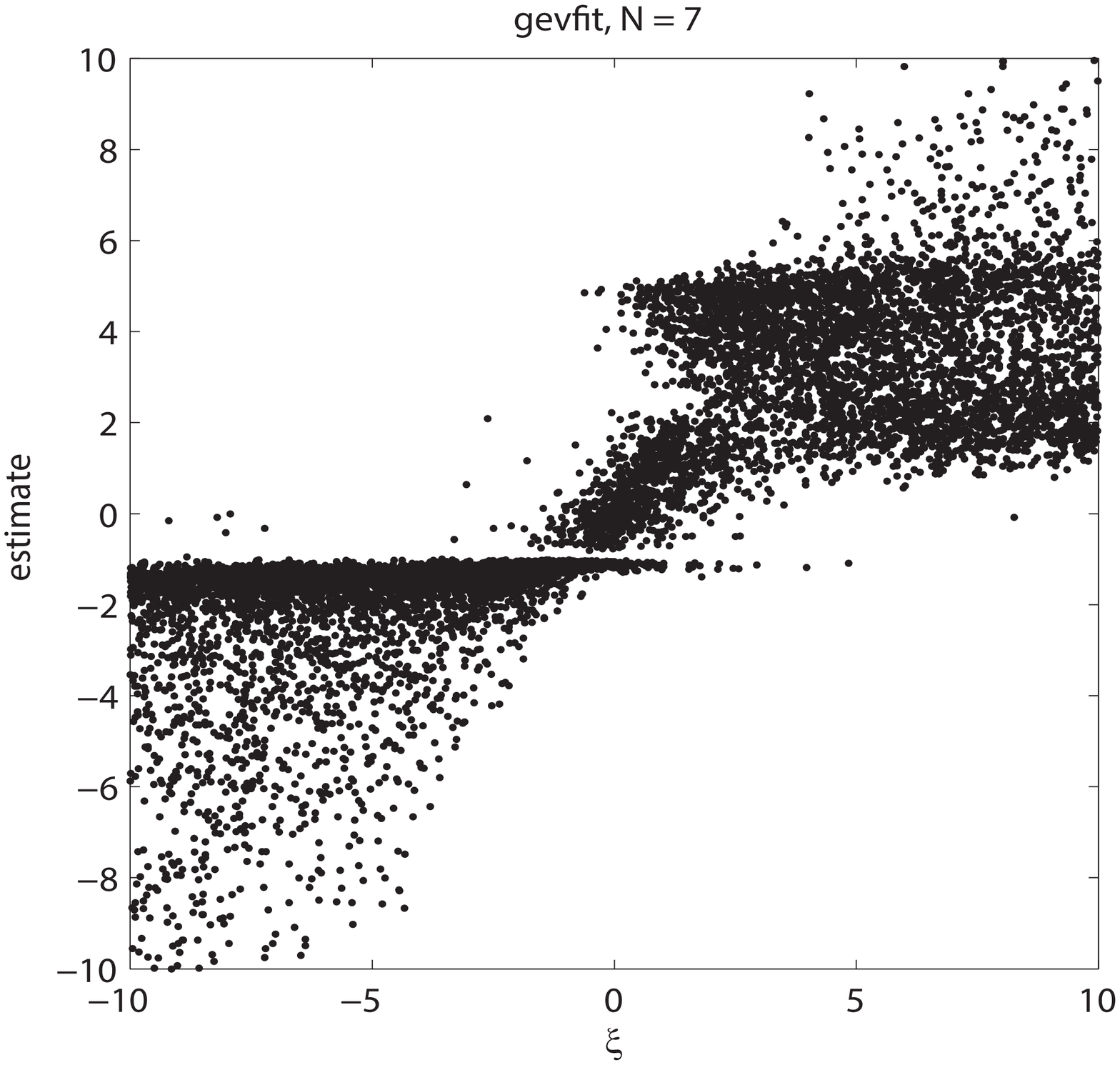}
   \includegraphics[width=70mm,keepaspectratio]{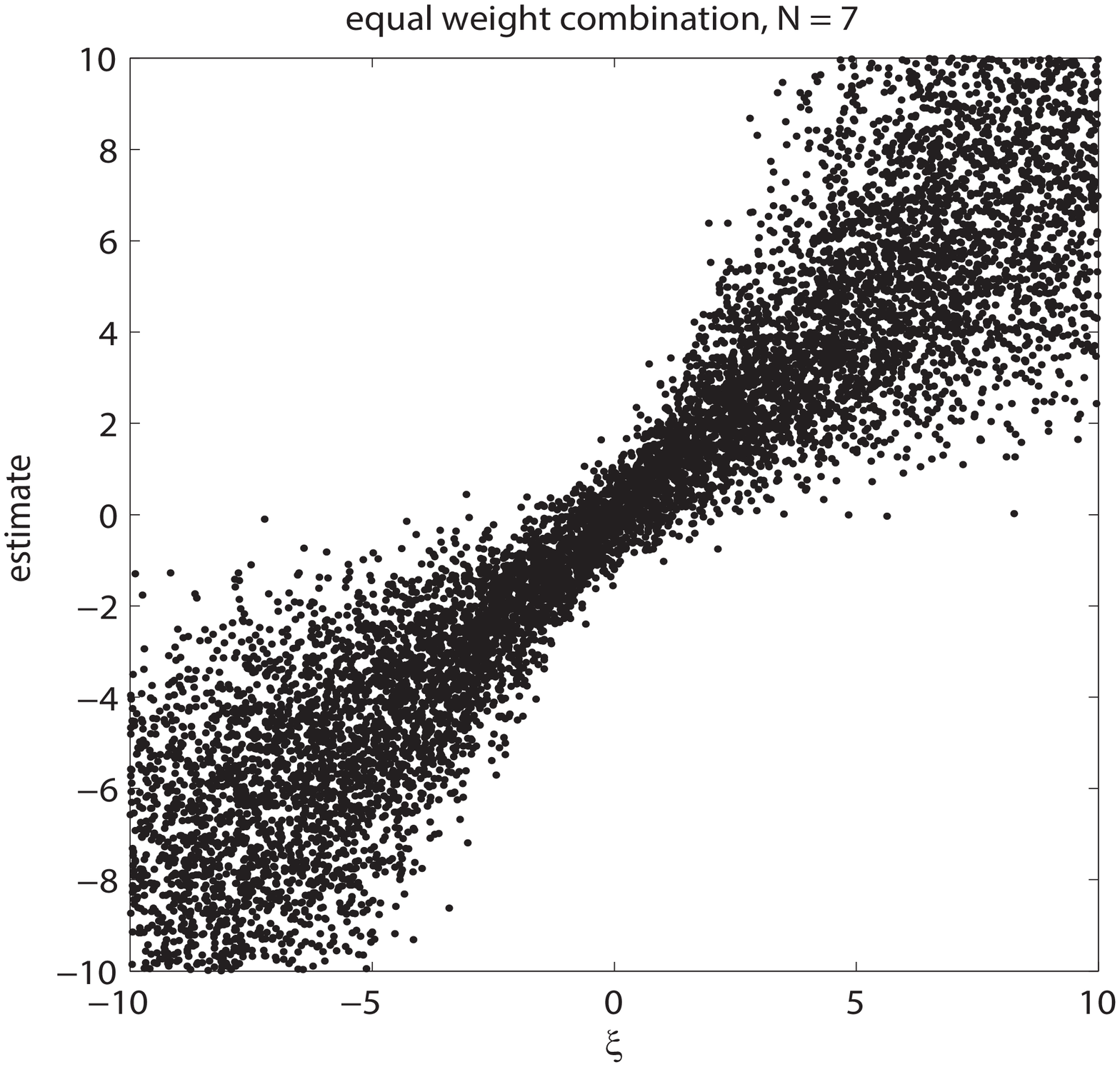}
  \caption{Each diagram shows 10,000 samples of size $N=7$, with $(\mu,\sigma) = (0,1)$ and $\xi$ uniform over
  $[-10,10]$. To the left are the resulting estimates obtained using Matlab {\it gevfit}, and to the right using the equal weight combination of elementals.
 }\label{GEVN7MLML}
\end{figure}

Fig.~\ref{GEVN7MLML} shows the estimates obtained via Maximum Likelihood and via the equally-weighted combination of elementals
for 10,000 samples of size $N=7$ drawn from a GEV. The parameters are $(\mu,\sigma) = (0,1)$ and $\xi$ is uniform random over
  $[-10,10]$. The Maximum Likelihood algorithm used is the {\it gevfit} function within the Statistics Toolbox of Matlab (\cite{matlab}).
The iterative ML analysis took considerably longer to run than the straight-forward elemental evaluation, and the unusual pattern of the ML
results suggests that the algorithm experienced some numerical difficulties. The elemental results, on the other hand, are as one may expect - showing unstructured scatter about a diagonal mean.

%

\begin{figure}[h!] \centering
  \includegraphics[width=120mm,keepaspectratio]{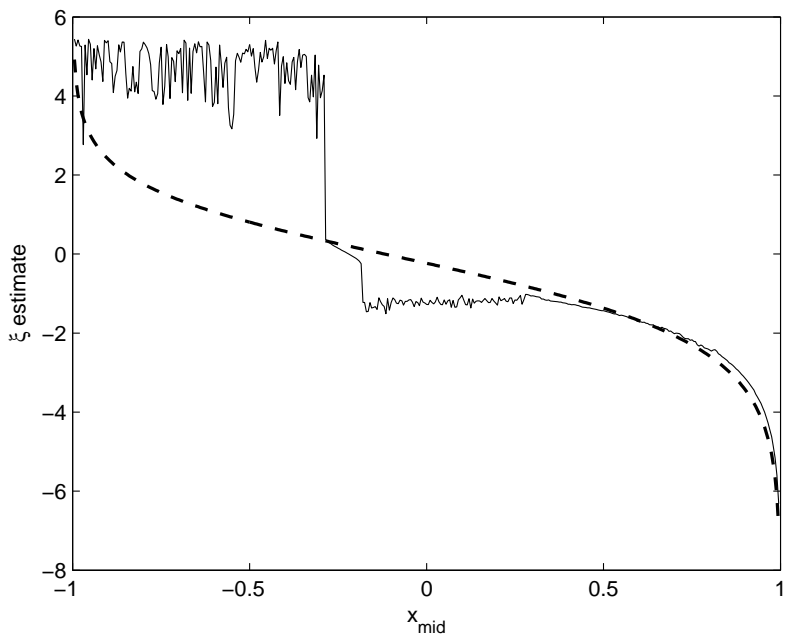}
  \caption{The estimates obtained using Matlab's {\it gevfit} ML estimator (solid) and the unique elemental (dashed) for the sample $[-1,x_{mid}, 1]$, as $x_{mid}$ is varied.
 }\label{MLmini}
\end{figure}

Statements about the efficacy of ML which are predicated on the parameters lying within some range - such as the regularity condition requiring $\xi < -1/2$ - are somewhat
unhelpful, given that the analyst typically does not know the parameters, but has only data.
As an example, consider an ordered  data set consisting of three points $[-1,x_{mid},1]$ where $x_{mid}$ is some value between -1 and 1.
Such data could have originated from a GEV of any $\xi$, thus questions such as whether $\xi < -1/2$ are not well posed when given only the data.
As $x_{mid}$ is varied, the estimates of $\xi$ using the {\it gevfit} ML algorithm and the (unique) elemental combination are shown in Fig. \ref{MLmini}.
Again, the irregular graph of the ML estimates is indicative of numerical difficulties, and there were indeed numerous error messages warning of convergence to boundary points or of
failure to converge within some
pre-designated number of iterations. By contrast, the elemental estimator produces a smooth, well-behaved graph whose computation requires no iteration, being a simple one-step
evaluation of the function
$\hat{\xi} = a \log \tau  -b \log t$ with $a = 1/(3\log 4/3 )$  and $b = 1/(3\log3/2)$. Moreover, the resulting graph of the elemental estimate conforms reasonably with
intuition, in that the centrally-located  $x_{mid} = 0$ leads to the estimate $\xi = -(\log^2 2 )/(3\log(3/4)\log(2/3)) = -0.23$, which is not very far from the uniform distribution
case of $\xi = -1$. (Indeed, there is no reason that the estimate for this uniformly-spaced data case should deliver the exact value $\hat{\xi} = -1$, but some loose
agreement is to be expected.)
Values of $x_{mid}$ to the right of this correspond to increasingly negative estimates of $\xi$, which correspond to the data tending to cluster towards the rightmost end point.
Similarly, values of $x_{mid}$ further left of centre lead to increasingly positive values of $\xi$, in correspondence with the tendency of positive $\xi$ distributions to cluster the data to the left end-point.
If one assumes that the elemental estimate is in some sense sensible, then it is perhaps surprising to observe that the ML estimate agrees closely with this for large negative values of $\xi$, the very parameter regime that regularity considerations would have led us to be wary of.

As further demonstration of the robustness of the elemental approach as compared with the ML approach, the estimator performance is compared over a number of idealised data sets.
To create an idealised sample of size $N$, the unit interval is divided into $N$ equal segments, the midpoints $F_i$  of which are
used, via the Probability Integral Transform, to create a sample $\lbrace x_i \ | \ x_i = F^{-1}(F_i), \ i = 1, \ldots, N  \rbrace$ where $F$ is a GEV distribution function with parameters
$(\mu, \sigma, \xi)$ with
$\mu = 0$, $\sigma = 1$, and $\xi$ varying over the interval $[-10,10]$. (The ordering of $i$ is such that $x_1$ is the data maximum).
Although each data set employs a nominal value of $\xi$ in its construction,
any such data set could have arisen by random sampling from a GEV with any other value of $\xi$.
The value of $\xi$ used in the idealised construction is thus only nominally associated with that data set.

\begin{figure}[h!] \centering
  \includegraphics[width=70mm,keepaspectratio]{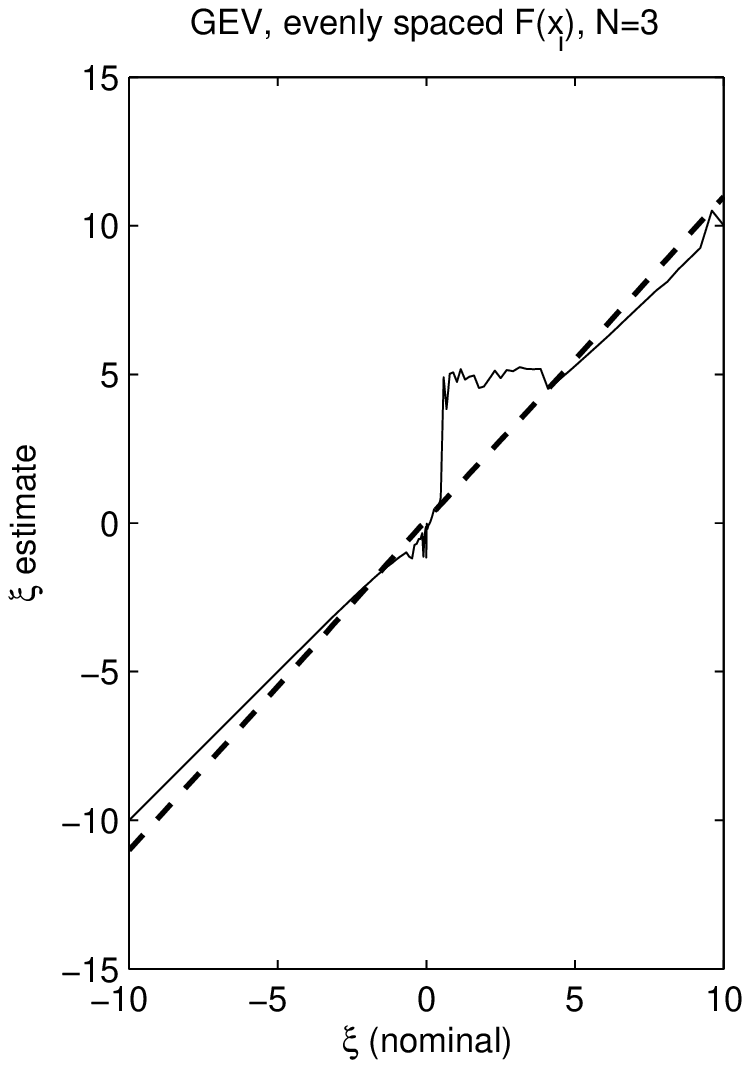}
  \includegraphics[width=70mm,keepaspectratio]{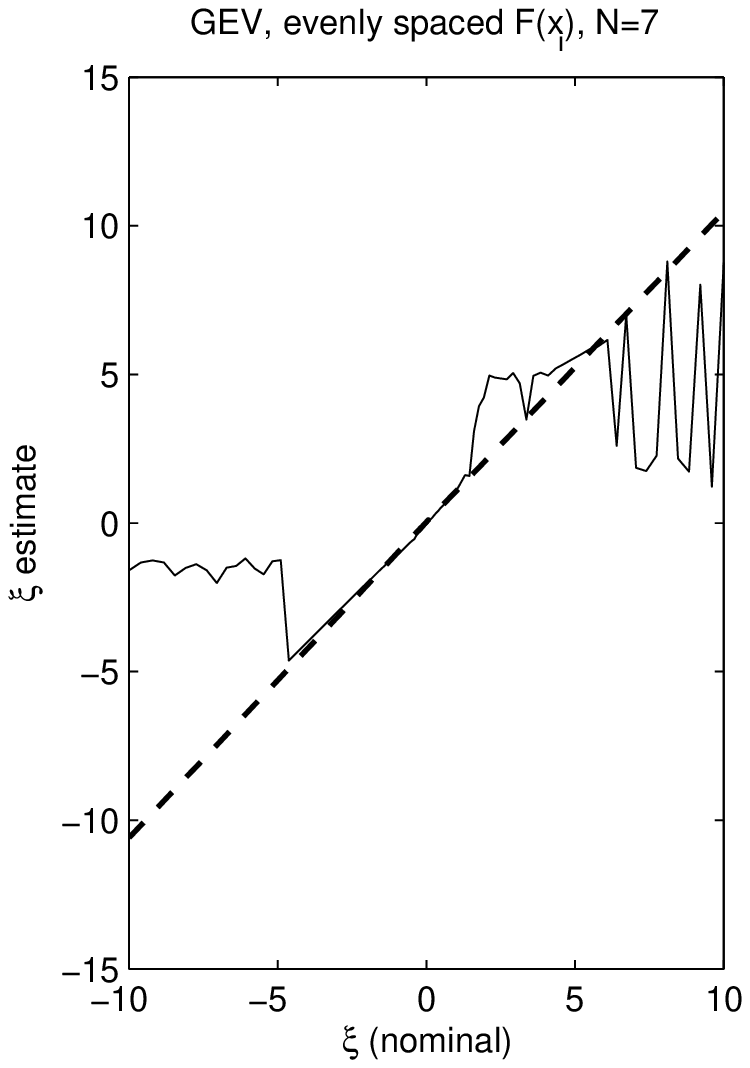}\\
  \includegraphics[width=70mm,keepaspectratio]{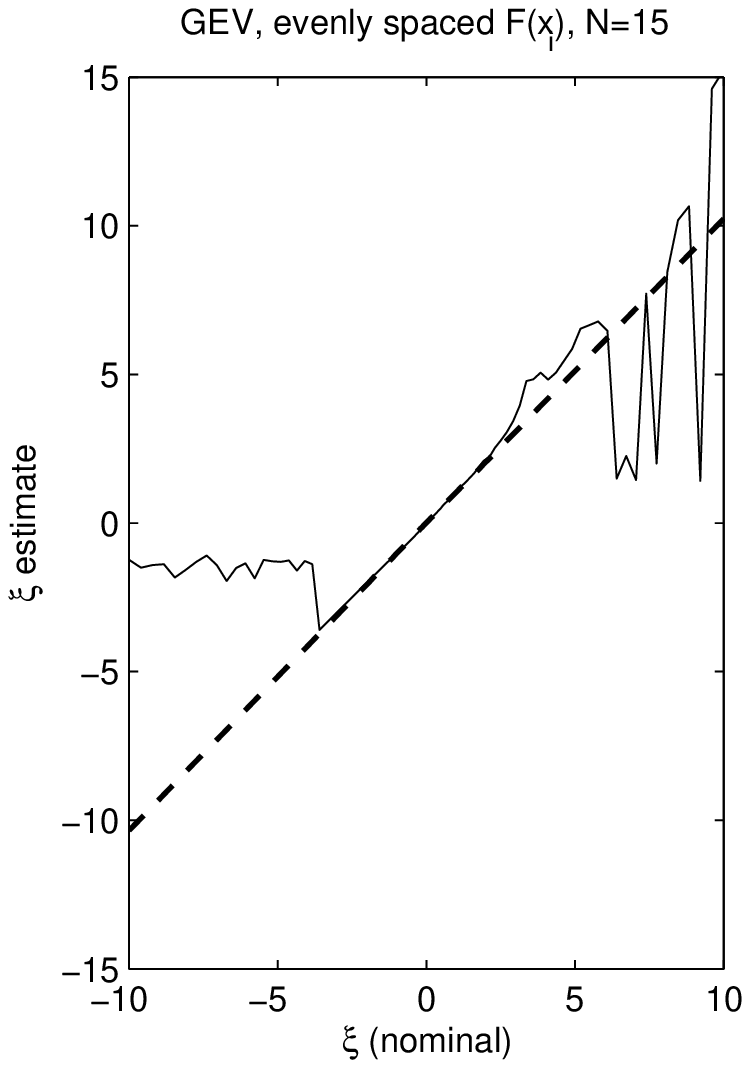}
  \includegraphics[width=70mm,keepaspectratio]{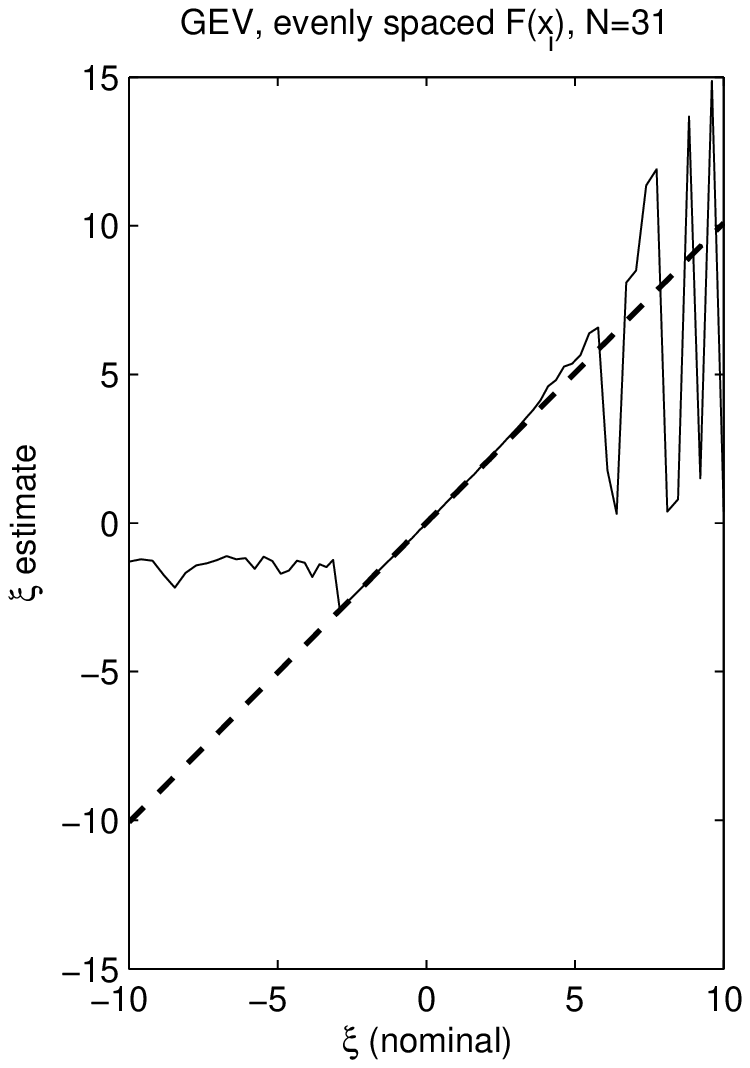}
  \caption{The estimates obtained using the Matlab {\it gevfit} ML estimator (solid) and the equal weight combination of  elementals (dashed) for the idealised samples created by uniformly-spaced $F$, for $N= 3,7, 15$ and 31.
 }\label{MLmini7N3}
\end{figure}

Fig.~\ref{MLmini7N3} compares the resulting {\it gevfit} ML estimates with those obtained from the equal weight elemental combination, for $N = 3,7, 15$ and 31.  Again, the irregularities of the ML graphs
are indicative of numerical difficulties. That these irregularities are present even for $N =31$ shows that they are not merely a feature of small sample sizes.
The elemental estimates however, are well-behaved for all samples, and lie close to the diagonal in all cases. Indeed, given the artificiality of the idealised data construction,
although some approximate diagonal correspondence is to be expected,
perfect diagonality is not - especially for $N$ small.

In summary, linear combinations of elemental estimators provide computationally efficient and robust estimators for the tail parameter,
and - even for extreme parameter values ($|\xi| \approx 10$) and very small sample sizes ($N \approx 3$) -  do not suffer the numerical difficulties experienced by the ML estimator.

\section{Application to distributions related to Weibull distributions}
If the density function of a GEV is flipped right to left, we obtain a three-parameter distribution related to the Weibull distribution.
It follows immediately that the shape parameter of a sample drawn from such a distribution
can be estimated by simply reversing the elemental coefficients.
That is, given a sample (ordered such that $X_1$ is  maximum) drawn from the three parameter $(\mu,\sigma,\zeta)$ distribution with distribution function
 \begin{equation}
F(x) = \left\{  \begin{array}{ll}
 1- \exp \left( - \left( 1 - \zeta z \right)^{-1/\zeta} \right) & \ \ \text{for} \ \ \zeta \neq 0 \\
 1 - \exp \left( - \exp \left(  z  \right) \right) & \ \ \text{for} \ \ \zeta = 0
\end{array} \right. \ \ \text{where} \ \ z = \frac{x-\mu}{\sigma}
\label{Weibull}
\end{equation}
 an elemental estimator of the parameter $\zeta$ can be obtained via
\begin{equation}
\hat{\zeta}_{IJN} = a_{N}^W(J) \log \tau   - b_{N}^W(I) \log t   \label{eqnweib}
\end{equation}
with $\tau$ and $t$ defined as per Eqn.~\ref{taut} and
\begin{equation}
a_{N}^W(J) = -b_{N}(N+1-J) \ \ \ \ \text{and} \ \ \ \  b_{N}^W(I) = -a_{N}(N+1-I)
\end{equation}
For large $N$, $b_N^W(I)$ converges to $ I \log (I/N)$.

The practical significance of all this is that we are often interested in heavy-tailed distributions which head off unbounded to the right, such as the GEV with positive $\xi$. By flipping the GEV right to left, we obtain another family of heavy-right-tailed distributions (and these are related to the Weibull distribution). These heavy right tails are those that, before flipping, headed off to the left in the GEV for $\xi$ negative. These tails have a different structure to the positive $\xi$ GEV tails, and we now have estimators for the shape parameter of this new class of heavy tails. However the shape parameter $\zeta$ of these new tails should not be confused with their extreme value tail index, the tail parameter $\xi$ of the GEV within whose domain of attraction this Weibull-related distribution lies. That tail index is 0, since the Weibull-related distribution lies within the domain of attraction of the Gumbel distribution, the GEV with $\xi = 0$ (\cite{embrechts}).

A similar reflection could be applied to the GPD estimators of \cite{McRobieGPD} but this is of arguably lesser interest, given that all left tails in the GPD are bounded below, and thus when flipped, lead to right tails that are bounded above.

\section{Summary}
Elemental estimators have been presented for the GEV tail parameter $\xi$. These have small bias, are computationally simple and robust, are applicable to small data sets and are valid for all parameters $\xi$, positive, zero and negative. Linear combinations of elementals appear to provide efficient, consistent estimators when applied to data drawn from pure GEVs, and such combinations appear to have advantages over Maximum Likelihood approaches to the GEV. Further, by a simple reflection, the approach has been shown to be applicable to another class of heavy-tailed distributions related to the Weibull distribution.



\bibliography{McRobieGEVbib}

\section*{Appendix 1. Derivation of the elemental coefficients}

The Generalised Extreme Value distribution has distribution function
\begin{equation}
F(x) = \left\{ \begin{array}{cc}
                 \exp \left[ - (1+ \xi \frac{x-\mu}{\sigma})^{-1/\xi} \right] & \ \ \text{for} \ \xi \neq 0 \\
                 \exp \left[ - \exp(-(\frac{x- \mu}{\sigma})) \right]  & \ \  \text{for} \   \xi = 0
               \end{array} \right.
\end{equation}
The support is
$\mu-\sigma/\xi \leq x$ for $\xi$ positive and $  x \leq \mu - \sigma / \xi$ for $\xi$ negative.
For now, we ignore the $\xi = 0$ case.
Inverting the distribution function via the Probability Integral Transform gives
\begin{equation}
x = u(F) = \mu + \frac{\sigma}{\xi}\left[ (-\log F)^{-\xi} -1 \right]
\end{equation}

Let $\bX$ be an ordered sample of size $N$,  ($X_N \leq X_{N-1} \leq \ldots
\leq X_2  \leq X_1$). The expected value of any function $h(\bX)$ is then
\begin{equation}
\langle h(\bX) \rangle =  N! \ \int_0^1 dF_1 \ldots \int_0^{F_{N-1}}
dF_N \ h(u(\bF))
\end{equation}
the integral being over the $N$-dimensional unit simplex defined by the ordering.

For any function of just two data points $X_i$ and $X_j$,  $(i < j)$, all other variables can be integrated out to leave
\begin{equation}
\langle h(X_i,X_j) \rangle = c_{ijN} \int_0^1 \ dF_i \int_0^{F_i}
dF_j \  \ \ f_{ijN}(F_i, F_j) \ h(u(F_i), u(F_j))
\end{equation}
where
\begin{equation}
c_{ijN} = \frac{N!}{(i-1)! (j-i-1)!(N-j)! }  \ \ \ \ \ \text{and} \ \ \ \ f_{ijN}(x,y) = (1-x)^{i-1} (x-y)^{j-i-1} y^{N-j}
\label{CIJN}
\end{equation}

The expected value of any function $h(X_k)$ of a single data point, is likewise, after a further integration, given by
 \begin{equation}
\langle h(X_k) \rangle = d_{kN} \int_0^1 \ dF_k \ \ g_{kn}(F_k) \ h(u(F_k))
\end{equation}
where
\begin{equation}
d_{kN} = \frac{N!}{(k-1)! (N-k)! }  \ \ \ \ \ \text{and} \ \ \ \ g_{kN}(x) = x^{N-k} (1-x)^{k-1}
\label{dkN}
\end{equation}

Expressed in terms of $F$'s, a log-spacing is given by
\begin{equation}
\log(X_i - X_j) = \left\{ \begin{array}{cc}
                            \log(\frac{\sigma}{\gamma}) + \gamma \log(-\log F_j) + \log \left[ 1- \left( \frac{-\log F_i}{-\log F_j} \right)^\gamma \right]  &  \text{for} \ \xi = -\gamma , \gamma > 0\\
                            \log(\frac{\sigma}{\gamma}) - \gamma \log(-\log F_i) + \log \left[ 1- \left( \frac{-\log F_i}{-\log F_j} \right)^\gamma \right]  &   \text{for} \ \xi = +\gamma , \gamma > 0
                          \end{array} \right.
\end{equation}

Already, the strategy is becoming evident. The zero sum of log-spacing weights in each elemental will remove the $\log (\sigma/\gamma)$ terms, and expectation of the $ \log(-\log F)$ terms will provide some constant multiple of $\gamma$ that will form the core of the estimator. One hopes that, as in the derivation for the GPD elementals,  the remaining complicated terms $ \log \left[ 1- \left( (-\log F_i)/(-\log F_j) \right)^\gamma \right]$ can be combined in such a way as to vanish.
Unfortunately this is not the case, although - as will be shown - the residual error is small.

\vspace*{5mm}

We shall seek an estimator of the form
\begin{equation}
\hat{\xi}_{IJN} = A \log \tau - B \log t
\end{equation}
with $\tau$ and $t$ defined as in Eqn.~\ref{taut}.

For $\xi$ negative ($\xi = -\gamma$, $\gamma$ positive), we obtain almost immediately
\begin{eqnarray}
\langle \log t \rangle_- & =  & \langle \log \left( \frac{X_{I+1} - X_J}{X_I - X_J} \right) \rangle \nonumber \\
                         & =  &  \langle \log \left[ 1 - \left( \frac{-\log F_{I+1}}{-\log F_J} \right)^\gamma \right] \rangle -
                                 \langle \log \left[ 1 - \left( \frac{-\log F_{I}}{-\log F_J} \right)^\gamma \right] \rangle \nonumber \\
                         & = &  \int_0^1 dx \int_0^x dy  \left[ c_{I+1,JN}f_{I+1,JN}(x,y) -   c_{IJN}f_{IJN}(x,y) \right] \log \left( 1- \left( \frac{-\log x}{- \log y}  \right)^\gamma \right) \nonumber \\
                         & \equiv  & I_1
\end{eqnarray}

Essentially there are only the complicated terms, denoted as $I_1$.

Similarly, continuing with $\xi$ negative, we obtain
\begin{eqnarray}
\langle \log \tau \rangle_- & =  & \langle \log \left( \frac{X_{I} - X_{J-1}}{X_I - X_J} \right) \rangle \nonumber \\
                         & =  &  \gamma \langle \log (-\log F_{J-1}) \rangle  - \gamma  \langle \log (-\log F_{J}) \rangle   \nonumber \\
                          & &  \ \ \ \    + \langle \log \left[ 1 - \left( \frac{-\log F_{I}}{-\log F_{J-1}} \right)^\gamma \right] \rangle -
                                 \langle \log \left[ 1 - \left( \frac{-\log F_I}{-\log F_{J}} \right)^\gamma \right] \rangle \nonumber
\end{eqnarray}
As earlier, the complicated terms can be collected into a single integral, by defining
\begin{equation}
                I_2 \equiv   \int_0^1 dx \int_0^x dy  \left[ c_{I,J-1,N}f_{I,J-1,N}(x,y) -   c_{IJN}f_{IJN}(x,y) \right] \log \left( 1- \left( \frac{-\log x}{- \log y}  \right)^\gamma \right)
\label{I2defn}
\end{equation}
The leading $\log(-\log F)$ terms are each a function of a single data point, and thus via Eqn.~\ref{dkN}
\begin{equation}
\langle \log \tau \rangle_- = \gamma \int_0^1 \ \left[d_{J-1,N} g_{J-1,N}(x) - d_{JN} g_{JN}(x) \right] \log( - \log x) \ dx    + I_2
\end{equation}

Denoting the term in square brackets as $\left[ DG(x) \right]$, we have
\begin{eqnarray}
\left[  DG  \right] & = & \frac{N!}{(J-2)!(N-J+1)!} x^{N-J+1} (1-x)^{J-2}  - \frac{N!}{(J-1)!(N-J)!} x^{N-J} (1-x)^{J-1}  \nonumber \\
                          & = &  \frac{N!}{(J-1)!(N-J+1)!} \left[ (J-1) x^{N-J+1} (1-x)^{J-2}  - (N-J+1)  x^{N-J} (1-x)^{J-1} \right]  \nonumber
\end{eqnarray}

Each of the $(1-x)^p$ terms can be expanded using the Binomial Theorem to obtain
\begin{equation}
\left[  DG  \right] = \left(
                              \begin{array}{c}
                                N \\
                                J-1 \\
                              \end{array}
                            \right)
\left\{
\sum_{m=0}^{J-2} \frac{(J-1) (J-2)!}{m! (J-2-m)!}(-1)^m x^{N-J+1+m} -
\sum_{m=0}^{J-1} \frac{(N-J+1) (J-1)!}{m! (J-1-m)!}(-1)^m x^{N-J+m}
\right\}
\end{equation}
In the first summation, we relabel the dummy summation index $m$ using $ k = m+1$ , and then relabel $k$ as $m$  (i.e. $m$ becomes $m-1$) giving
\begin{equation}
\left[  DG \right] = \left(
                              \begin{array}{c}
                                N \\
                                J-1 \\
                              \end{array}
                            \right)
\left\{
\sum_{m=1}^{J-1} \frac{(J-1) (J-2)!}{(m-1)! (J-1-m)!}(-1)^{m-1} x^{N-J+m} -
\sum_{m=0}^{J-1} \frac{(N-J+1) (J-1)!}{m! (J-1-m)!}(-1)^m x^{N-J+m}
\right\} \nonumber
\end{equation}
Collecting terms for each $m =1$ to $J-1$, we obtain
\begin{equation}
\left[  DG \right] = \left(
                              \begin{array}{c}
                                N \\
                                J-1 \\
                              \end{array}
                            \right)
\left\{
-(N-J+1) x^{N-J} -
\sum_{m=1}^{J-1} (-1)^m x^{N-J+m} \frac{(J-1)!}{m! (J-1-m)!} \left(  N-J+1+m \right)
\right\}
\end{equation}
Noting that the leading $-(N-J+1) x^{N-J}$ term is that which would be obtained if the summation were extended down to $m = 0$, we obtain
\begin{equation}
\left[ DG \right] = - \left(
                              \begin{array}{c}
                                N \\
                                J-1 \\
                              \end{array}
                            \right)
\left\{
\sum_{m=0}^{J-1} (-1)^m x^{N-J+m}
\left(
                              \begin{array}{c}
                                J-1\\
                                m \\
                              \end{array}
                            \right)
 \left(  N-J+1+m \right)
\right\}
\end{equation}

We now need only multiply by $\log(-\log x)$ and integrate over the unit interval, and use the result that
\begin{equation}
\int_0^1 x^{\alpha -1} \log (-\log x ) \ dx = \frac{-\log \alpha}{\alpha} - \frac{C}{\alpha}
\end{equation}
where $C$ is Euler's constant.

The terms involving Euler's constant are an alternating sum of binomial coefficients and thus sum to zero:
\begin{equation}
\sum_{m=0}^{J-1} (-1)^m
\left(
                              \begin{array}{c}
                                J-1\\
                                m \\
                              \end{array}
                            \right)
\ \ = 0
\end{equation}

This leaves only the $- (\log \alpha)/\alpha$ terms (where $\alpha = N-J+1+m$), which give
\begin{equation}
\int_0^1    \left[ DG(x) \right] \log( -\log x ) \ dx
= \left(
                              \begin{array}{c}
                                N \\
                                J-1 \\
                              \end{array}
                            \right)
\left\{
\sum_{m=0}^{J-1} (-1)^m
\left(
                              \begin{array}{c}
                                J-1\\
                                m \\
                              \end{array}
                            \right)
\log( N-J+1+m)
\right\}
 \equiv \frac{-1}{A}
\end{equation}

Collecting all the above, we obtain finally that for $\xi$ negative,
\begin{equation}
\langle \log \tau \rangle_- = \frac{-\gamma}{A} + I_2 = \frac{\xi}{A} + I_2
\end{equation}

We now need to repeat the procedure for $\xi$ positive. Writing $ \xi = \gamma$, $\gamma >0$, we readily obtain
that
\begin{equation}
\tau = \frac{1 - \left( \frac{-\log F_I}{-\log F_{J-1}}\right)^\gamma}{1 - \left( \frac{-\log F_I}{-\log F_{J}}\right)^\gamma}
\end{equation}
whose expected logarithm was defined in Eqn.~\ref{I2defn} as $I_2$, thus
\begin{equation}
\langle \log \tau \rangle_+  = I_2
\end{equation}

Similarly for $\langle \log t \rangle_+$, we obtain an $I_1$ contribution from the complicated integrals, and the $\log (-\log F)$ terms are
\begin{equation}
\gamma \left[ \langle \log(-\log F_I) \rangle - \langle \log(-\log F_{I+1}) \rangle  \right]
\end{equation}
These are identical to the terms in the $\langle \log \tau \rangle_-$ derivation, involving integrations over a single variable, except we must make the substitution $J-1$
goes to $I$.
We thus obtain
\begin{equation}
\langle \log t \rangle_+  = \frac{-\xi}{B} + I_1
\end{equation}
with
\begin{equation}
\frac{1}{B} =
 \int_0^1 \left[ \ \ldots \ \right] \log ( -\log ( x) \ dx  = - \left(
                              \begin{array}{c}
                                N \\
                                I \\
                              \end{array}
                            \right)
\left\{
\sum_{m=0}^{I} (-1)^m
\left(
                              \begin{array}{c}
                                I\\
                                m \\
                              \end{array}
                            \right)
\log( N-I+m)
\right\}
\end{equation}

In summary,
\begin{eqnarray}
\langle \log t \rangle_-  & = &    \ \ \ \ \ \ \   I_1 \nonumber \\
\langle \log \tau \rangle_- & = & \frac{\xi}{A} + I_2 \nonumber \\
\langle \log t \rangle_+  & = & \frac{-\xi}{B} + I_1 \nonumber \\
\langle \log \tau \rangle_+  & = &  \ \ \ \ \ \ \    I_2 \nonumber \\
\end{eqnarray}
whence, for positive or negative  $\xi$, we have
\begin{equation}
A \langle \log \tau \rangle - B \langle \log t \rangle = \xi + A I_2 - B I_1
\end{equation}
This is an unbiased estimator for $\xi$ provided that the weighted sum of the complicated terms $A I_2 - B I_1 = 0$.
Although the equivalent terms in the GPD derivation (\cite{McRobieGPD}) cancelled exactly, this is not the case for the GEV here.
However the resulting bias is small. It also vanishes as $|\xi|$ becomes large.

\pagebreak

\section*{Appendix 2. Approximation of the $b_N(I)$ coefficients}
\begin{eqnarray}
\text{Writing} \ \ \frac{b_N(I)}{N} & = & \left[ N \left( \begin{array}{c}
                                         N \\
                                         I
                                       \end{array}  \right) S_I \right]^{-1} \ \ \ \ \text{with} \ \ \
S_I =  - \sum_{m=0}^I \left(\begin{array}{c}
                              I \\
                              m
                            \end{array} \right)
(-1)^m \log(N-I+m) \\
\text{and} \ \ \ \ f(x) & = & -(1-x)\log(1-x)  \ \ \ \ \text{with} \ \ \ x = \frac{I}{N}
\end{eqnarray}
we show that
\begin{equation}
\frac{b_N(I)}{N}  = f(x) \left[ 1 + \frac{x}{12N} + O\left(\frac{1}{N^3}\right) \right]
\end{equation}

The demonstration begins by writing the ratio of the approximation $Nf(x)$ to the true value of the coefficient
$b_N(I)$ as
\begin{equation}
\frac{f(x)}{b_N(I)/N} = \left[ N \left( \begin{array}{c}
                                         N \\
                                         I
                                       \end{array} \right) \right] S_I f(x)
\end{equation}
and we expand each of the three terms on the right-hand side as a power series in $1/N$, keeping the first three nonzero terms in each expansion.

\subsubsection*{Term 1:}
\begin{eqnarray}
N \left( \begin{array}{c}
                                         N \\
                                         I
                                       \end{array} \right)
                                       & = &  \frac{N.N(N-1)(N-2) \ldots (N-I+1)}{I!} \
                                       \ = \ \ \frac{N^{I+1}}{I!} \left(1-\frac{1}{N}\right)\left(1-\frac{2}{N}\right) \ldots \left(1-\frac{I-1}{N}\right) \nonumber \\
& = &  \frac{N^{I+1}}{I!} \left[
1 - \frac{I(I-1)}{2N} + \frac{I(I-1)(I-2)(3I-1)}{24N^2} + \orderNcube \right]
\end{eqnarray}

\subsubsection*{Term 2:}
\begin{eqnarray}
S_I & = &   - \sum_{m=0}^I \left(\begin{array}{c}
                              I \\
                              m
                            \end{array} \right)
(-1)^m \log(N-I+m)   \ \ \ =  \ \  \  (-1)^{I+1} \sum_{k=0}^I \left(\begin{array}{c}
                              I \\
                              k
                            \end{array} \right)
(-1)^k \log\left(1- \frac{k}{N} \right)    \nonumber \\
 & = &   (-1)^{I} \sum_{k=0}^I \left(\begin{array}{c}
                              I \\
                              k
                            \end{array} \right)
(-1)^k \sum_{j=1}^\infty \frac{1}{j}  \left( \frac{k}{N} \right)^j  \ \ \ \ = \ \ \
   (-1)^I  \sum_{j=1}^\infty \frac{1}{jN^j}
       \sum_{k=0}^I (-1)^k \left(\begin{array}{c}
                              I \\
                              k
                            \end{array} \right)
k^j
\end{eqnarray}

We now use the relations
\begin{equation}
(-1)^I \sum_{k=0}^I (-1)^k \left(\begin{array}{c}
                              I \\
                              k
                            \end{array} \right)
k^j  =      \begin{cases}
                   \ \ \ 0 \hspace*{45mm}  \text{for}  \ \ \ \ j = 1, \ldots, I-1   \\
                   \ \ \  I! \hspace*{44mm}   \text{for} \ \ \ \ j = I \\
                   \ \left[ I(I+1)/2 \right] \ I!   \hspace*{28mm}  \text{for} \ \ \ \ j = I+1  \\
                   \ \left[I(I+1)(I+2)(3I+1)/24 \right] \ I!   \hspace*{6mm} \text{for} \ \ \ \ j = I+2
                 \end{cases}
\end{equation}
giving
\begin{eqnarray}
S_I & = & \frac{I!}{IN^I} + \frac{I(I+1)I!}{2(I+1)N^{I+1}} + \frac{I(I+1)(I+2)(3I+1)I!}{24(I+2)N^{I+2}} +
{\large{O}}\left( \frac{1}{N^{I+3}}\right)  \nonumber \\
 & = &  \frac{I!}{IN^I}  \left[ 1 + \frac{I^2}{2N} + \frac{I^2(I+1)(3I+1)}{24N^2} +
{\large{O}}\left( \frac{1}{N^3}\right) \right]
\end{eqnarray}

This result is worthy of comment. Essentially we have expanded each $\log(1-k/N)$ as a power series in $k/N$. When these series are weighted by
the binomial coefficients and summed then - rather remarkably - the first $I-1$ terms disappear.
For example, if calculating $b_N(I)$ with $I = 100$, the first 99 terms in these series disappear, and
the first non-zero terms are those involving $(k/N)^{100}$.

\subsubsection*{Term 3:}
\begin{eqnarray}
f(x) = -(1-x)\log(1-x) &  = &  -(1-x)\left[ -x -\frac{x^2}{2} - \frac{x^3}{3} + {\large{O}}\left( x^4 \right) \right] \nonumber \\
& = & x -\frac{x^2}{2} - \frac{x^3}{6} + {\large{O}}\left( x^4 \right) \ \ = \ \ \frac{I}{N} \left[ 1- \frac{I}{2N} - \frac{I^2}{6N^2} + \orderNcube \right]
\end{eqnarray}

\subsubsection*{Combine terms}
Each of the three terms is now in the form $c_i (1 + a_i/N + b_i/N^2 + O(1/N^3))$ thus their product
can be written
 \begin{equation}
\frac{f(x)}{b_N(I)/N} = \left[ N \left( \begin{array}{c}
                                         N \\
                                         I
                                       \end{array} \right) \right] S_I f(x) = C \left[ 1 + \frac{A}{N} + \frac{B}{N^2} + \orderNcube \right]
\end{equation}
with
\begin{eqnarray}
C & = & c_1 c_2 c_3 = \left( \frac{N^{I+1}}{I!} \right) \left( \frac{I!}{IN^I} \right) \left( \frac{I}{N} \right) = 1 \nonumber \\
A  & = & a_1+a_2+a_3 = \frac{-I(I-1)}{2} +\frac{I^2}{2} +\frac{-I}{2} = 0             \nonumber \\
\text{and} \ \ \ B & = &   a_1a_2 + a_2a_3 + a_3 a_1 + b_1 + b_2 + b_3 \nonumber  \\
 & = &  \left[\frac{-I(I-1)}{2} \right]\left[\frac{I^2}{2} \right] +  \left[\frac{I^2}{2}\right]\left[\frac{-I}{2} \right] + \left[\frac{-I}{2} \right]\left[\frac{-I(I-1)}{2} \right] \nonumber  + \\
 & &  \hspace*{10mm} + \frac{I(I-1)(I-2)(3I-1)}{24} + \frac{I^2(I+1)(3I+1)}{24} - \frac{I^2}{6} \nonumber \\
 & = & \frac{1}{24} \left[ (-6I^4 +6I^3) +(-6I^3) +(6I^3-6I^2) + (3I^4-10I^3+9I^2 -2I) +(3I^4+4I^3+I^2) +(-4I^2)
 \right] \nonumber \\
  & = & \frac{-I}{12}
\end{eqnarray}
Finally, then
\begin{equation}
\frac{f(x)}{b_N(I)/N} = 1 - \frac{I}{12N^2}+ \orderNcube
\end{equation}
or
\begin{equation}
\frac{b_N(I)}{N} = -(1-x)\log(1-x) \left[ 1 + \frac{I}{12N^2}+\orderNcube \right]
\end{equation}

Thus, for $I$ fixed, the coefficient $b_N(I)$ converges to $ -N(1-x)\log(1-x)$ as $N$ increases.

The fore-going demonstration would be improved if the convergence were established at fixed $x$ rather than fixed  $N$. However,
we note, without further proof, that
\begin{equation}
\frac{b_N(I)}{N} \approx -(1-x)\log(1-x) - \frac{x}{12N} \log(1-x)
\label{better}
\end{equation}
is a better approximation and this is illustrated in Fig.~\ref{Acalc1a}.

The following examples illustrate the somewhat unusual nature of this approximation.
Consider the case $N=20$, $I = 3$.
\begin{equation}
b_{20}(3) = \frac{-3.2.1}{20.19.18} \frac{1}{\log \left( \frac{17.19^3}{18^3.20}\right)} \approx - \left( 17 + \frac{3}{12 (20)}\right) \log\left( \frac{17}{20}\right)
\end{equation}
for which a hand-held calculator gives $ 2.76494701 \approx 2.76494671$. The error is $\approx 3\times 10^{-7}$,  which is substantially better than the $O(1/N^3)$ suggested by the derivation.

Although the approximation was developed for large $N$, it works surprisingly well even for the smallest possible $N$ (namely $N=3$),
leading to somewhat unusual results such as
\begin{eqnarray}
b_3(1) = \frac{-1}{3\log(2/3)} & \approx & -\left( 2+ \frac{1}{12(3)}\right) \log(2/3)  \hspace*{10mm}
\text{or} \ \ \ \ 0.8221 \approx  0.8222   \nonumber \\
b_3(2)  = \frac{-1}{3\log\left(\frac{1 (3)}{2^2}\right)} & \approx & -\left(1+\frac{2}{12(3)}\right) \log(1/3)
\hspace*{10mm}
\text{or} \ \ \ \ 1.1587 \approx  1.1596   \nonumber
\end{eqnarray}

In summary, whilst it is conceptually simple to write down and evaluate the exact expression for $b_N(I)$ in terms of
 binomially-weighted sums of logarithms, the numerical instabilities are such that the
 remarkably accurate approximation of Eqn.~\ref{better} provides a useful alternative.

\end{document}